\documentclass[11pt]{amsart}

\usepackage{lineno,hyperref}
\usepackage{amsmath, amsthm, amssymb, enumerate}
\usepackage{color}
\usepackage{datetime}
\usepackage{fancyhdr}
\usepackage{amssymb,amsmath,amsthm,amsfonts,graphicx,hyperref,url,dsfont}

\usepackage{listings}
\usepackage{xcolor}
\usepackage[figure]{hypcap}

\usepackage{siunitx}
\usepackage{booktabs}
\usepackage[english]{babel}
\usepackage{blindtext}
\RequirePackage{fix-cm}

\newtheorem{thm}{Theorem}
\newtheorem{lem}[thm]{Lemma}
\newtheorem{rmk}{Remark}

\newtheorem{de}{Definition}[section]

\newtheorem{cor}{Corollary}[section]
\newtheorem{ex}{Example}[section]

\def\R{\mathbb R}
\def\C{\mathbb C}
\def \E{\mathcal E}
\def\N{\mathbb N}
\def \F{\mathcal F}
\def \P{\mathbb P}
\def \B{\mathcal B}

\def\D{\mathcal D}
\def\L{\mathcal L }
\def\cal{\mathcal }
\DeclareMathOperator*{\esssup}{ess\,sup}

\def\tilde{\widetilde}

\def\supp{\mathop{\rm supp}\nolimits}    
    
\usepackage[T1]{fontenc}
\usepackage[utf8]{inputenc}
\usepackage{combelow}

\usepackage{newunicodechar}

\usepackage{blindtext}
\usepackage{float}

\begin{document}
\title[Stochastic parabolic   equations with singular potentials]{Stochastic parabolic   equations with singular potentials}

\author[S. Gordi\'c, T. Levajkovi\'c, Lj. Oparnica]{Sne\v zana Gordi\'c$^1$, Tijana Levajkovi\'c$^2$, \\ Ljubica Oparnica$^3$}
\date{\today}

\thanks{$^1$Sne\v zana GORDI\' C, {Faculty of Education, University of Novi Sad, Podgori\v cka 4, 25000 Sombor, Serbia, \texttt{snezana.gordic@pef.uns.ac.rs}}\\
	\indent$^2$Tijana LEVAJKOVI\'C, {Institute of Stochastics and Business Mathematics, Vienna University of Technology,  Austria, \texttt{tijana.levajkovic@tuwien.ac.at}}\\
\indent$^3$Ljubica OPARNICA, Department of Mathematics: Analysis, Logic and Discrete Mathematics, University of Ghent, Krijgslaan 281 (Building S8), B 9000 Ghent, Belgium,  \texttt{oparnica.ljubica@ugent.be} \& Faculty of Education, University of Novi Sad, Podgori\v cka 4, 25000 Sombor, Serbia
}

\maketitle
\begin{abstract}
In this work we consider a class of  stochastic parabolic  equations with singular space depending potential, random driving force and random initial condition. For the analysis of these equations we combine the  chaos expansion method from the white noise analysis and the concept of very weak solutions. For given stochastic parabolic  equation we introduce the notion of a stochastic very weak solution, prove the existence and uniqueness of the very weak solution to corresponding stochastic  initial value problem and show its independence of  a regularization on  given singular potential. In addition, the consistency of a stochastic very weak solution with a stochastic  weak solution is shown. 
\end{abstract}

\section{Introduction} \label{sec1} 

We consider stochastic evolution problems with singular potentials which arise in probabilistic modelling of uncertainty in engineering and science, for example in structural mechanics, material science, fluid dynamics, climate and turbulence modelling. In particularly, the motivation comes from aggregation models in chemical kinetics, population dynamics, image processing, the modelling of options  in financial mathematics, pressure diffusion in a porous medium and aerodynamics. The aim of this work is  to analyse  stochastic evolution problems of the form
\begin{equation}\label{Eq:Mostgeneral}
\begin{split}
\left(\frac{\partial}{\partial t}  - {\mathcal L} \right) \, U  + \,q\, \cdot \, U &= F,\\
U |_{t=0} &= G,
\end{split}
\end{equation}
where $U$ is unknown, $F$ and $G$ are given generalized stochastic processes  depending on space, time and random component, $\L$ is an elliptic operator acting on the space variable only  (its action  on a stochastic process is interpreted  as an action on its space component).  The potential  $q$ is singular, and   depending on its nature, one  can interpret the product in \eqref{Eq:Mostgeneral} differently. We are interested in all possible singular behaviours of $q$, either in space and time or in the random component.  Typical examples  would be space and/or time white noise, inverse squared potential in space, Dirac delta distribution in space and/or time. The most general case is to consider  $q$ to be a random, space and time depending potential, i.e., a generalized  stochastic process, so that \eqref{Eq:Mostgeneral} takes the form
\begin{equation}
\label{Eq: general}
\begin{split}
\left(\frac{\partial}{\partial t}  - {\mathcal L} \right) \, U(t, x, \omega) + q( t, x, \omega) \, \lozenge \, U(t, x, \omega) &= F(t, x, \omega),\\
U(0, x, \omega) &= G(x, \omega).
\end{split}
\end{equation}
Here, $\lozenge$ denotes the Wick product and it is introduced  to give a sense to the product of two generalized stochastic processes: the potential $q$ and the unknown stochastic process $U$. In order to deal with singularities in deterministic variables  we employ the concept of  very weak solutions. 
The Wick product, also known as  stochastic convolution,  is the highest order stochastic approximation of the ordinary product. When at least one of the processes in the Wick product is a deterministic  function or an adapted process, it  reduces to the ordinary product. Alternative approaches for solving stochastic partial differential equations (SPDEs) with singularities have been developed in the theory of regularity structures  \cite{hairer2014} and in rough path theory and paracontrolled distributions \cite{gubinelli2004,GIP2015}. Another possibility is to consider the equation in Colombeau algebras of generalized functions and after regularization interpret the product as a classical product \cite{GOPS2019,OR-C2005,OR1998,Russo1994,RO1999,Selesi2008}. 

An example of \eqref{Eq: general} is the stochastic heat equation with random potential.  
Heat equations with random input data  are extensively studied in the literature due to their various applications in  biology, aerodynamics, structural acoustics, financial mathematics \cite{AHR1997,ANV2000,BDP1997,KL2017}. 
The heat equation with random potential, also known as the Anderson model, appears in the context of chemical kinetics and population dynamics \cite{hairer2015}.
In \cite{hu2019} a survey of recent progress on stochastic heat equations driven by Gaussian noise is given. Semilinear heat equation which is driven by a space-time Gaussian white noise in suitable algebras of generalized functions is considered in \cite{RO1999}. Stochastic evolution problems with polynomial nonlinearities were studied in \cite{LPSZ2018}. 
We aim to consider  more general classes of problems, i.e., to study stochastic parabolic evolution equations with singular space depending potentials, where the operator $\L$ is not necessarily the Laplace operator.

For the analysis of \eqref{Eq: general} we propose a new  method which combines the chaos expansion method  \cite{HOUZ1996,hu2015,LPSZ2018,tudor2013,xiu2002}  and the very weak solution concept  \cite{ARST2020b,ARST2020c,ARST2021,Gar20,GR15,MRT2019a,MRT2019b,RT17a,RT17b,RT2018,RY2020}. 
The chaos expansion method  is based on constructing the solution to the SPDE as a Fourier series in terms of a Hilbert space basis of orthogonal stochastic polynomials, with unknown coefficients being elements in an  appropriate space of deterministic functions.  As a result,  the initial problem  \eqref{Eq: general} is reduced to a system of deterministic parabolic  equations with singular potentials.
To deal with strong singularities in these deterministic equations we employ the concept of  very weak solutions  introduced in \cite{GR15}.  The idea is to model irregular objects in equations by approximating nets of regular  functions with moderate assymptotics. One obtains a net of  regularized  problems which can be treated in a usual distributional  way. As a result we obtain a net of solutions, which if moderate is called very weak solution.  
We apply this  to each  of the deterministic parabolic equations arising from the chaos expansion method. 
The obtained nets of very weak solutions are the coefficients of the unknown stochastic  process.  Summing them up and proving its convergence in an appropriate space of stochastic processes, one obtains the solution to the initial stochastic problem. 

In this work we concentrate on the simplest form of \eqref{Eq: general}  that involves  singularities. We analyze the case when $q$ is a deterministic singular space potential and $F$ and $G$ are singular stochastic processes of Kondratiev-type.  Since $q$ is not random, the Wick product in \eqref{Eq: general} is then the ordinary product. 
We particularly study the initial value problem for the stochastic parabolic equations  of the form 
\begin{equation}
\label{Eq: stochastic evolution}
\begin{split}
\left(\frac{\partial}{\partial t}  - \mathcal{L} \right) \, U(t, x, \omega) + q(x) \cdot U(t, x, \omega) &= F(t, x, \omega), \enspace \\
U(0, x, \omega) &= G(x, \omega),
\end{split}
\end{equation}
where $t\in (0,T]$, $x\in \mathbb R^d$, $ \omega\in \Omega$.  The operator $\mathcal L$  is unbounded and closed operator on  $L^2(\R^d)$ with dense domain, which generates $C_0$-semigroup on $L^2(\R^d)$. For example,  $-\mathcal L$  could be an elliptic differential operator of even order or the Laplace operator. 
The potential $q$ is a distribution in space. For example, it can be of the form $q=q_1+\delta,$ where $\delta$ is the Dirac delta distribution and $q_1$ is an $L^\infty$ function over $\mathbb{R}^d$.

The paper is organized as follows. In the sequel  we briefly introduce  notation
and basic concepts used in the following sections. In Section 2 we study initial value problem for the 
stochastic parabolic  equations with bounded potential, and  the deterministic evolution equations with singular space potential. Section 3 is devoted to stochastic parabolic  equation 
with singular space depending potential. The new method for its analysis is introduced. In Section 4 we illustrate our method on the stochastic heat equation with singular potential, we list several advantages of the introduced approach and indicate further extensions. 

\newpage
\subsection{Notations and basic concepts} \label{NotNot}
This paper is placed in the framework of white noise analysis and is built on the well known Wiener-It\^o chaos expansion theorem that relates the Gaussian measure and Hermite polynomials  \cite{HOUZ1996,xiu2002}.  Throughout the paper  $(\Omega, \mathcal F, \mathbb P)$ denotes the Gaussian white noise space and  $T>0$.  

Let $L^\infty(\R^d)$ be the space of essentially bounded measurable functions, \linebreak $C^k([0,T])$ be the space of $k$-times continuously differentiable functions on  $[0,T],$  $\mathcal{D}(\R^d):=\mathcal{C}_{0}^{\infty} (\mathbb{R}^d)$  the space of compactly supported smooth functions, $\mathcal{D}'(\mathbb{R}^d)$  the space of distributions, $\mathcal{S} (\mathbb{R}^d)$  the Schwartz space of rapidly decreasing functions, $\mathcal{S}' (\mathbb{R}^d)$  the Schwarz space  of tempered distributions, and $\mathcal{E}'(\mathbb{R}^d)$  the space of compactly supported distributions. Further, let $L^2(\R^d)$ be the space of square integrable functions over $\R^d$,  $H^{1}(\mathbb{R}^d)$ be the Sobolev space, a subset of functions $f$ in $L^2(\mathbb{R}^d)$ such that $f$ and its weak derivatives have a finite $L^2$-norm, $H_0^1 (\R^d)$ the closure of $\mathcal{D}(\R^d)$ in  $H^{1}(\mathbb{R}^d)$, $L^2([0,T]; L^2(\mathbb R^d))$  the Banach space of square integrable functions over $[0,T]$  with values in $L^2(\mathbb R^d)$ with the norm 
$$
\|f\|_{L^2([0,T]; L^2(\mathbb R^d))} =\left( \int_0^T \|f(t, \cdot)\|^2_{L^2(\mathbb R^d)} \, dt \right)^{1/2}
=  \left( \int_{\mathbb R^d }  \int_0^T   |f(t, x) |^2   dt dx \right)^{1/2},$$ 
$C([0,T];  L^2(\mathbb R^d))$  the Banach space of continuous functions over $[0,T]$  with values in $L^2(\mathbb R^d)$ with the norm  $$\|f\|_{C([0,T]; L^2(\mathbb R^d))} = \sup_{t \in [0,T]} \|f(t, \cdot)\|^2_{L^2(\mathbb R^d)},$$
and $L^1(0,T; {L^{2}(\mathbb R^d)})$ is the space of integrable functions over $[0,T]$ with values in  $L^{2}(\mathbb R^d).$

The Gaussian white noise probability space  $(\Omega,\F,\P)$ is  $(\mathcal{S}'(\R^d),\B,\mu)$, where  $\B$ is the family of all Borel subsets of $\mathcal{S}'(\R^d)$ equipped with the weak$^*$ topology and $\mu$ is the Gaussian white noise probability measure, whose existence is guarantied by the Bochner--Minlos theorem \cite{HOUZ1996}. The  space of square integrable random variables on  $(\Omega,\F,\P)$, denoted by $L^2(\mu)$,   is a Hilbert space with the norm $\|\cdot\|_{L^2(\mu)}$ induced by the inner product 
\[(F, G)_{L^2(\mu)} = \mathbb{E} \, (FG) = \int_{\Omega} F(\omega) \, G(\omega) \, d\mu(\omega) <\infty, \quad \text{for}\,\,  F, G\in L^2(\mu),\]where $\mathbb{E}$ denotes the expectation with respect to the Gaussian measure $\mu$. 

Let $\mathcal{I} :=\N^m_0 $ be the set of sequences of non-negative integers which have finitely many nonzero components $\gamma= (\gamma_1,\gamma_2, \ldots,\gamma_m,0,0,\ldots), \gamma_i \in \N_0, i=1, 2,\ldots, m, m\in\N,$ the $k$th unit vector $e(k)= (0,\ldots, 0,1,0, \ldots), k\in \N$  is the sequence of zeros with the number $1$ as the $k$th component and  $\textbf{0}= (0,0,\ldots)$  the zero vector,  $|\gamma|= \sum_{i=1}^{\infty}\gamma_i$ the length of multi-index $\gamma$ and $\gamma ! = \prod_{i=1}^{\infty} \gamma_i !$.

For a given $\gamma\in \mathcal{I}$, the $\gamma$th {\it Fourier-Hermite polynomial} is defined by
\begin{equation}\label{Halpha}
H_{\gamma}(\omega)=\prod_{k=1}^\infty h_{\gamma_k}(\langle\omega,
\xi_k\rangle), \quad\gamma\in\mathcal I, 
\end{equation}
where $\xi_k$, $k\in \mathbb{N}$ is the Hermite function of order $k$  and $h_k$, $k \in \mathbb N_0$ are the Hermite polynomials. 
The famous  Wiener-It\^o chaos expansion theorem states that each square integrable random variable  $F\in L^2(\mu)$ has a unique representation of the form
\begin{equation*}
F(\omega)  = \sum_{\gamma \in \mathcal{I}} \, c_{\gamma} \, H_{\gamma}(\omega),
\, c_\gamma\in \mathbb{R}, \omega\in \Omega, \quad \text{with}\quad
\|F\|^2_{L^2(\mu)}=\sum_{ \gamma\in \mathcal{I}} \, c_{\gamma}^2 \, \gamma! \, <\infty . 
\label{L2_rv_cond} 
\end{equation*}

The space of Kondratiev test random variables $(S)_{p}, p\in \N_0$ consists of the elements 
$f (\omega)= \sum_{\gamma \in \mathcal{I}} c_\gamma H_\gamma(\omega) \in L^2(\mu),$ $c_\gamma \in \R,$ such that 
$	\Vert f \Vert_{p}^{2} = \sum_{\gamma \in \mathcal{I}}^{} c_{\gamma}^{2}(\gamma !)^{2} (2\N)^{\gamma p} < \infty, $ where 	$(2\N)^\gamma = \prod_{i=1}^{\infty} (2i)^{\gamma_i},$ $\gamma \in \mathcal{I}$. 
The Kondratiev space of random variables $(S)_1$ is the projective limit of the spaces $(S)_{ p}$, i.e.,  $(S)_1 = \bigcap_{p \in \N_0}(S)_{ p}$. The family of seminorms $\Vert \cdot \Vert _{ p}$ generates a topology on $(S)_1$. 
The Kondratiev space of  stochastic distributions $(S)_{-p} $, consists of the formal expansions 
$F(\omega)=\sum_{\gamma \in \mathcal{I}} b_\gamma H_\gamma (\omega),$ $b_\gamma \in \R$ such that
\begin{equation*}
\Vert F\Vert_{-p}^{2} = \sum_{\gamma \in \mathcal{I}} b_{\gamma}^{2} (2\N)^{-\gamma p} < \infty.  \label{kondr seminorms gen}
\end{equation*}
The Kondratiev space $(S)_{-1}$ is the inductive  limit of the spaces $(S)_{-p}$, i.e.,  $(S)_{-1} = \bigcup_{p \in \N_0}(S)_{ -p}$. It is the dual space of $(S)_1$ and   a Frech\'et space.  
The action of $F\in(S)_{-1}$ on $f\in(S)_1$ is given by
$
\langle F, f \rangle = \sum_{\gamma \in \mathcal{I}}  b_\gamma  c_\gamma \gamma !  \in \R.
$
The spaces  $(S)_{ p}$ and $(S)_{ -p}$ are separable  Hilbert spaces, $(S)_1$ and $(S)_{-1}$ are nuclear, i.e.,  the embedding $(S)_{q} \subseteq (S)_{p}$ for $p\leq q$ is Hilbert-Schmidt, and
for $p\leq q$ it holds 
$
(S)_{ 1} \subseteq (S)_{ q} \subseteq (S)_{ p}  \subseteq L^2(\mu) \subseteq   (S)_{ -p}  \subseteq   (S)_{ - q} \subseteq   (S)_{ - 1}. 
$

For any normed space $Y,$  the tensor product $Y\otimes  L^2(\mu)$ is the space of  $Y$- valued square integrable stochastic processes, $Y\otimes (S)_{ -p }$ the space of $Y$-valued generalized stochastic processes  of Kondratiev-type and the space $Y\otimes (S)_{-1}$ is the inductive limit of  spaces $Y\otimes (S)_{-p }$, $p\geq 0$. The Wiener--It\^o chaos expansion theorem  can be extended  from spaces of random variables to spaces of  stochastic processes, \cite{LM2017, LPSZ2015, LS2017,PS2007}. 
We summarize known results in the following theorem.  
\begin{thm}Let  $Y$ be a normed space. 
	\begin{itemize}	
		\item[(a)] It holds
		\begin{equation}
		\label{2N konv p>1}
		\sum_{\gamma\in \mathcal I} (2\N)^{-p\gamma} < \infty \qquad \text{ if and only if}  \qquad p>1. 
		\end{equation}   
		\item[(b)] Every $Y$-valued square integrable stochastic process $F\in Y\otimes  L^2(\mu)$ can be uniquely represented in the chaos expansion form 
		\begin{equation}\label{CE}
		F(\omega) = \sum_{ \gamma\in \mathcal I}f_\gamma  \, H_\gamma(\omega),\quad  f_\gamma\in Y,
		\end{equation}where
		$\|F\|^2_{Y\otimes L^2(\mu)} = \sum_{ \gamma\in \mathcal{I}}\|f_\gamma\|^2_Y \,  \gamma!<\infty$. 
		If $Y$ is a Banach space, then $Y\otimes L^2(\mu)$ is a Banach space.
		\item[(c)] Every $Y$-valued generalized stochastic process  of Kondratiev-type $F \in Y\otimes (S)_{-1}$ can be  represented in the chaos expansion form \eqref{CE}  where
		\begin{equation*}\label{ConvergenceCondition}
		\|F\|^2_{Y\otimes (S)_{-p }} = \sum_{ \gamma\in \mathcal{I}}\|f_\gamma\|^2_Y  \, (2\N)^{-p\gamma}<\infty 
		\end{equation*}
		holds for some $p\geq 0$. 
		If $Y$ is a Banach space, then $Y\otimes(S)_{-p}$ are  Banach spaces for every $p\geq 0$, and $Y\otimes (S)_{-1}$ is a Frech\'et space.
	\end{itemize}
\end{thm}
Let $Y$ be a space of deterministic functions depending on time and space and  $U=U(t, x, \omega)\in Y\otimes (S)_{-1}$.  For fixed $\omega \in \Omega$ the process  $U(\cdot, \cdot, \omega)$ is a deterministic function in $Y$, and for fixed  $t\in(0,T]$ and $x\in\mathbb R^d$ the process $U(t, x, \cdot)$   belongs to a Kondratiev space of stochastic random variables.

\begin{ex}\label{Example0}
	\begin{enumerate}
		\item[(i)]  An $C^k([0, T]) $-valued  generalized process   $U\in C^k([0, T]) \otimes (S)_{-1}$ has the chaos expansion representation \[U(t, \omega) = \sum_{\gamma\in \mathcal I} u_\gamma(t) \, H_\gamma(\omega), \quad t\in [0,T], \, \omega\in \Omega,\] with the  coefficients $u_\gamma$, $\gamma\in \mathcal I$,  being elements of    the Banach space of functions $C^k([0, T])$, $ k\in \mathbb N$ such that for some $p \geq 0$
		\[\|U\|^2_{C^k([0, T])\otimes (S)_{ -p }} = \sum_{ \gamma\in \mathcal I}\|u_\gamma\|^2_{C^k([0, T])}  \, (2\N)^{-p\gamma}<\infty.\]
		\item[(ii)]  Let $\mathcal{B}$ be a Banach space of functions depending on $t$ and $x$. An $C([0,T]; \mathcal{B}) $-valued Kondratiev-type  generalized stochastic process  $U\in C([0,T]; \mathcal{B})\otimes (S)_{-1} $ has a chaos expansion representation 
		\begin{equation}\label{ChaosExpansion}
		U(t, x, \omega) = \sum_{\gamma\in \mathcal I} u_\gamma(t, x) \, H_\gamma(\omega), \quad t\in [0,T], \, x\in \mathbb R^d, \, \omega\in \Omega,
		\end{equation}
		with $u_\gamma\in C([0,T]; \mathcal{B}),$
		$\gamma\in \mathcal I$, and 	for some  $p\geq 0$ it holds 
		\begin{equation*} 
		\|U\|^2_{C([0,T];\mathcal{B}) \otimes (S)_{ -p}} = 
		\sum_{ \gamma\in \mathcal I}\|u_\gamma\|^2_{C([0,T]; \mathcal{B})}  \, (2\N)^{-p\gamma}	< \infty.
		\end{equation*}
		\item[(iii)] An $L^2([0,T], L^2(\mathbb R^d))$-valued generalized process   $U\in L^2([0,T], L^2(\mathbb R^d))$ $ \otimes (S)_{-1}$ has the chaos expansion representation 
		\eqref{ChaosExpansion}
		with the coeffi\-cients  $u_\gamma\in L^2([0,T], L^2(\mathbb R^d))$, $\gamma\in \mathcal I$, and 	for some  $p\geq 0$ it holds 
		\[\begin{split} &\|U\|^2_{L^2([0,T], L^2(\mathbb R^d)) \otimes (S)_{ -p}} = 
		\sum_{ \gamma\in \mathcal I}\|u_\gamma\|^2_{L^2([0,T], L^2(\mathbb R^d))}  \, (2\N)^{-p\gamma}	< \infty.
		\end{split}\]
		\item[(iv)] The time white noise $W$ is  an element of  the space $C^{k}([0,T])\otimes (S)_{-1}$ and it is given formally by	
		\begin{equation}
		\label{whitenoise}
		W(t,\omega)=\sum_{k=1}^{\infty} \xi_k(t)H_{e(k)}(\omega),
		\end{equation} 
		while space-time white noise belongs to   $C^{k}([0,T],\mathbb{R}^d) \otimes (S)_{-1}$ and 
		\begin{equation*}
		W(t,x,\omega)=\sum_{k,n=1}^{\infty} \xi_k(t) \eta_n (x) H_{e(k)}(\omega),
		\end{equation*}
		where $\xi_k,$ $k\in \mathbb{N}$ are the Hermite functions and $\eta_n$, $n \in \mathbb{N}$ are elements of an orthogonal basis of  $\mathbb{R}^d$,  see \cite{HOUZ1996}.
	\end{enumerate}
\end{ex}

\section{Set up of the problem and methodology}
We study the  stochastic parabolic  initial value problem \eqref{Eq: stochastic evolution},  i.e., 
\begin{equation*}
\begin{split}
\left(\frac{\partial}{\partial t}  - \mathcal{L} \right) \, U(t, x, \omega) + q(x) \cdot U(t, x, \omega) &= F(t, x, \omega), \enspace \\
U(0, x, \omega) &= G(x, \omega),
\end{split}
\end{equation*}
where $t\in (0,T]$, $x\in \mathbb R^d$, $ \omega\in \Omega$. 
For $X\subseteq L^2(\R^d)$ we define ${\cal X}(X): =  C([0,T]; X) \otimes (S)_{-1}$. Note that if $X$ is a Banach space then ${\cal X}(X)$ is Banach space as well. Assume the following:
\begin{enumerate}
	\item[(H1)] The operator $ \mathcal L $ is unbounded and closed operator on  $L^2(\R^d)$ 	with dense domain $D\subseteq L^2(\mathbb R^d)$,  which generates  a $C_0$-semigroup    $(T_t)_{t\geq 0}$ on $L^2(\R^d)$.
	Action of the operator $ \mathcal L $  on a generalized  stochastic  process  of Kondratiev-type $U\in {\cal X}(D)$ with the chaos expansion  \eqref{ChaosExpansion} is given by
	\[{\mathcal L} U(t, x, \omega) := \sum\limits_{\gamma\in \mathcal I} \L \, u_\gamma(t, x) \, H_\gamma(\omega),\]
	where $\L$ acts only on the space component.
	\item[(H2)] The force term $F\in {\cal X}(L^2(\R^d))$  is $C([0,T]; L^2(\R^d)) $-valued Kondratiev-type  generalized stochastic process with the chaos expansion \eqref{ChaosExpansion} where the coefficients $f_\gamma\in C([0,T]; L^2(\R^d))$, $\gamma\in \mathcal I$	 are Lipschitz continuous functions with respect to $t$.
	\item[(H3)] The initial condition $G$  is a $D$-valued Kondratiev-type  generalized stochastic process, i.e.,  $G\in D \otimes (S)_{-1}$. 
\end{enumerate}

\vspace{0.2cm}
\subsection{Stochastic parabolic  equation with bounded potential} \label{subsec 2.2}
We start our analysis of  (\ref{Eq: stochastic evolution}) assuming $q\in L^\infty(\R^d)$.
The following theorem employs the chaos expansion method in order to show existence of the unique solution to (\ref{Eq: stochastic evolution}) for $q\in L^\infty(\R^d)$.
\begin{thm}
	\label{stochastic weak sol} 
	Let the operator $\mathcal L$, the force term $F$ and the initial condition $G$ satisfy the assumptions $(H1)$-$(H3)$.  Assume the potential $q \in L^\infty(\mathbb R^d)$. Then,  there exists unique  generalized stochastic process  $U \in {\cal X}(D)\subseteq {\cal X} ({L}^2(\R^d)) $ satisfying stochastic parabolic initial value problem  (\ref{Eq: stochastic evolution}).
\end{thm}

{\bf Proof.}
	By assumptions all stochastic processes appearing in \eqref{Eq: stochastic evolution} are of Kondra\-tiev-type, and can be  represented in their chaos expansions  \eqref{ChaosExpansion}. For $t\in [0,T]$,  $x\in \mathbb R^d$, and  $\omega\in \Omega$ we have
	\begin{equation}
	\label{F chaos}
	F(t, x, \omega)=\sum\limits_{\gamma\in \mathcal{I}}f_\gamma(t, x)  H_\gamma(\omega),
	\end{equation}
	with 
	\begin{equation}
	\label{force term uslov konv}
	\sum\limits_{\gamma\in \mathcal I} \|f_\gamma\|_{C([0,T]; L^2(\R^d)) }^2  \,\, (2\mathbb N)^{-p_1\gamma}< \infty
	\end{equation}for some  $p_1\geq 0$
	and
	\begin{equation}
	\label{G chaos}
	G(x, \omega) = \sum\limits_{\gamma\in \mathcal I} g_\gamma(x) \, H_\gamma(\omega), 
	\end{equation} with
	\begin{equation}\label{poc uslov konv}
	\sum\limits_{\gamma\in \mathcal I} \|g_\gamma\|_{L^2(\mathbb R^d)}^2  \, \, (2\mathbb N)^{-p_2\gamma}< \infty
	\end{equation}
	for some $p_2\geq 0$. 
	Assume that the solution is also given in the  chaos expansion form \eqref{ChaosExpansion}, i.e.,
	$$
	U(t, x, \omega) = \sum\limits_{\gamma\in \mathcal I} u_\gamma(t, x) \, H_\gamma(\omega),\quad t\in [0,T],\,x\in \mathbb R^d,\,\omega\in \Omega.
	$$
	The chaos expansion method means to substitute stochastic processes given in their chaos expansion forms into the equation and to equalize the corresponding coefficients with respect to orthogonal stochastic polynomial basis $H_\gamma(\omega)$, $\gamma\in \mathcal I$. The initial stochastic problem \eqref{Eq: stochastic evolution} is then reduced to a system of the deterministic PDEs   of the form  
	\begin{equation}
	\label{problemm Pgammaa}
	\begin{split}
	\left(\frac{\partial}{\partial t} -  {\L}  \right)u_\gamma(t, x) + q(x) \cdot u_\gamma(t, x) &= f_\gamma(x, t),\\
	u_\gamma(0, x) &= g_\gamma(x),
	\end{split}
	\end{equation} for every $\gamma \in \mathcal I$.  
	By hypotheses we have that   $q\in L^\infty(\mathbb R^d)$,    $f_\gamma \in C([0,T];L^2(\R^d))$  are Lipschitz continuous functions with respect to $t$ and $g_\gamma \in D$ for all $\gamma\in \mathcal I$. By Theorem  \ref{theorem estimates L infinity deterministic} (proved below independently of this theorem) it follows that for each $\gamma\in \mathcal I$ the  deterministic problem \eqref{problemm Pgammaa} has a unique bounded nonegative solution $u_\gamma \in C([0,T];D)\subseteq C([0,T];L^2(\R^d))$ given by
	\begin{equation}
	\label{semigroup solution u det 1r}
	u_\gamma (t,x)=   S_t g_\gamma (x)+\int_0^t S_{t-s} f_\gamma(s,x)\,ds, \quad t\in (0,T], \, x \in \mathbb{R}^d,
	\end{equation}
	and such that  for $t\in (0,T]$  the following estimate holds
	\begin{equation} 
	\label{ee-perturbationsr}
	\|u_\gamma(t,\cdot)\|_{L^2(\mathbb R^d)} \leq M(t) \left( \|g_\gamma(\cdot)\|_{L^{2}(\mathbb R^d)} +\int_0^t \|f_\gamma(s,\cdot)\|_{L^2(\mathbb R^d)} \,ds\right).
	\end{equation}
	Here $M(t)= M(t, w, M, \|q\|_{L^\infty(\mathbb R^d)}) = Me^{(w + M \|q\|_{L^\infty(\mathbb R^d)}) t} $ is an increasing positive function on $[0,T]$, with stability constants $w$ and $M$  appearing in the semigroup estimate \eqref{semigroup estimate} in the proof of Theorem \ref{theorem estimates L infinity deterministic}.
	
	Finally, we are going to show that the process \eqref{ChaosExpansion}
	whose  coefficients $u_\gamma$ are given by  \eqref{semigroup solution u det 1r} converges in ${\cal X}(D)$.
	Therefore we have to show that for some $r\geq 0$
	\begin{equation*}
	\|U\|^2_{X \otimes (S)_{-r}} = \sum_{\gamma\in \mathcal I} \|u_\gamma\|^2_X \, \,  (2\mathbb N)^{-r\gamma} < \infty, 
	\label{kondr type converg} 
	\end{equation*}
	where here (and along this proof) 
	$X:= C([0,T];L^2(\R^d))$. Namely, for $r\geq \max{\{p_1, p_2\}}$ by the estimates  \eqref{ee-perturbationsr} we obtain  
	\[\begin{split}
	\|U\|^2_{X \otimes (S)_{-r}} &\leq  2M(T)^2 \sum_{\gamma\in \mathcal I} \left( \|g_\gamma\|^2_{L^2(\mathbb R^d)} + 
	\left(\int_0^T \|f_\gamma(t, \cdot)\|_{L^2(\mathbb R^d)} dt\ \right)^2\right) (2\mathbb N)^{-r\gamma}\\
	&  \leq 2M(T)^2  \left(\sum_{\gamma\in \mathcal I}  \|g_\gamma\|^2_{L^2(\mathbb R^d)}  (2\mathbb N)^{-p_2\gamma}  + T^2 \sum_{\gamma\in \mathcal I}  \|f_\gamma\|^2_X \, (2\mathbb N)^{-p_1\gamma} \right)\\
	& \leq 2M(T)^2  \left(\|G\|^2_{L^2(\mathbb R^d)\otimes (S)_{ -p_2}} + T^2 \|F\|^2_{X \otimes (S)_{ -p_1}}\right),
	\end{split}\]
	which is by the assumptions $(H2)$ and $(H3)$ finite. 
	
	The uniqueness of the solution $U$ follows from the uniqueness of its coefficients $u_\gamma$, $\gamma\in \mathcal I$, which are given by \eqref{semigroup solution u det 1r} and the uniqueness of  the chaos expansion representation in the Fourier--Hermite basis of orthogonal stochastic polynomials. \hfill $\square$ \\

For the analysis of the problem \eqref{Eq: stochastic evolution} with {$q\in \D'(\R^d)$}  
in Section \ref{sec4} we will use the same approach. There the deterministic equation \eqref{problemm Pgammaa} for every $\gamma \in \mathcal{I}$ will have irregular potential. 

\vspace{0.3cm}
\subsection{A deterministic parabolic  equations with singular  potentials: The very weak solution approach} \label{sec3}

Let us now consider the  deterministic parabolic  equation of the form  \eqref{problemm Pgammaa}, with the potential  {$q\in \D'(\R^d)$},  i.e.,
\begin{equation}
\label{deterministic problem}
\begin{split}
\left( \frac{\partial}{\partial t}  -  {\cal L} \right)  u(t, x) + {q (x)} \cdot  u(t, x) & = f(t, x ), \quad t\in (0,T], \,   x\in  \mathbb R^d,  
\\
u(0, x) &= g(x),
\end{split}
\end{equation} 
where $\cal L$ is an operator on $L^2(\R^d)$ with dense domain $D\subseteq L^2(\R^d)$, as in the hypothesis (H1), $f\in C([0,T]; L^2(\R^d))$ Lipschitz continuous function with respect to $t$, and $g\in D$. 

\subsubsection{Bounded potential and related estimates}

First we consider the case with $L^\infty$-potential $q$ and provide an  estimate for the solution to \eqref{deterministic problem}. In particular, the following theorem assures that \eqref{problemm Pgammaa} in Theorem \ref{stochastic weak sol} has a unique bounded nonegative solution $u_\gamma \in C([0,T]; D) $, $\gamma\in\mathcal{I}$.

\begin{thm}
	\label{theorem estimates L infinity deterministic}
	Let the operator $\L$ be as in (H1), the force term $f\in C([0,T]; L^2(\R^d))$ a Lipschitz continuous function with respect to $t$, the initial condition $g\in D$,  and let the potential $q \in L^{\infty} (\mathbb{R}^d)$. Then,  the deterministic parabolic   initial value 
	problem \eqref{deterministic problem} has a unique bounded nonnegative solution  $u\in C([0,T]; D) )\subseteq C([0,T]; L^2 (\mathbb{R}^d))$ 
	satisfying   
	\begin{equation} 
	\label{ee-perturbations}
	\|u(t,\cdot)\|_{L^2(\mathbb R^d)} \leq  M(t)  \left( \|g(\cdot)\|_{L^{2}(\mathbb R^d)} +\int_0^t \|f(s,\cdot)\|_{L^2(\mathbb R^d)} \,ds\right)
	\end{equation}
	for $t\in (0,T],$	where 
	\begin{equation}
	\label{eModt}
	M(t) : = M \exp {\left( \left(w + M \|q\|_{L^\infty(\R^d)}\right)t \right)}, \quad t\in (0,T],
	\end{equation}
	with $w\in\R$ and $M>0$ being the stability constants from the semigroup estimates \eqref{semigroup estimate}.
\end{thm}

{\bf Proof.} 
	As the operator $\L$ is assumed to be the infinitesimal generator  of a $C_0$-semigroup 	we aim in applying the semigroup theory.  We rewrite the parabolic problem \eqref{deterministic problem} in the form 
	\begin{equation}
	\label{det sistem pregrupisan}
	\begin{split}
	\frac{\partial}{\partial t}  u(t, x) &=  \left( {\L} - {q (x)} \text{Id}\, \right) u(t, x) + f(t, x ), \quad t\in (0,T], \,   x\in  \mathbb R^d, 
	\\
	u(0, x) &= g(x),
	\end{split}
	\end{equation}
	with $\text{Id}$ denoting the identity operator on $L^2(\mathbb R^d)$.	 
	The operator $q \, \text{Id}$, i.e., {operator of multiplication of an element from $L^2(\mathbb R^d)$ by $q \in L^{\infty}(\R^d)$}, is bounded operator on $L^2(\mathbb R^d)$ with the bound $Q=\esssup_{x \in \R^d} |q(x)| $. Indeed, 
	\[\begin{split}
	\|q(\cdot) u(t,\cdot)\|_{L^2(\mathbb R^d)} & = \left(\int_{\R^d} |q(x)u(t,x)|^2\,dx\right)^{1/2} \\  &\leq   \esssup_{x \in \R^d} |q(x)|\left( \int_{\R^d} |u(t,x)|^2\,dx \right)^{1/2} \\
	& = Q \,  \|u(t,\cdot)\|_{L^2(\mathbb R^d)}.
	\end{split}\]
	The operator $\L$ is an infinitesimal generator of a $C_0$-semigroup of operators 	$(T_t)_{t\geq 0}$ on $L^2(\mathbb R^d)$ such that 
	\begin{equation}
	\label{semigroup estimate}
	\|T_t\|_{L(L^2(\R^d))}\leq Me^{wt}, \quad t\geq 0
	\end{equation}
	holds for some $M>0$ and $w\in\mathbb R$.
	Therefore, since $L^2(\R^d)$ is a reflexive Banach space, by 	\cite[Chapter 3, Theorem 1.1] {Pazy1983}
	the operator $\L - q \, \text{Id}$ is the infinitesimal generator of a $C_0$-semigroup $(S_t)_{t \geq 0}$ on $L^2(\mathbb R^d)$ satisfying 
	\begin{equation} 
	\label{semigroupp estimate 2}
	\|S_t\|_{L(L^2(\R^d)) }\, \leq \,  M \exp{\left(w+M\|q\|_{L^\infty(\R^d)}\right) t}, \quad t\in [0,T],
	\end{equation} 
	where $M$ and $w$ are as in \eqref{semigroup estimate}. 
	
	Since 
	$f\in C([0,T]; L^2(\R^d))\subset L^1(0,T;L^2(\R^d))$ is Lipschitz continuous function with respect to $t$, by the result from \cite[Chapter 4, Corollary 2.11]{Pazy1983}  for every $g\in D$   there exists a unique $u\in C([0,T]; D) $ which is differentiable almost everywhere on $[0,T]$   
	solving the deterministic parabolic problem \eqref{det sistem pregrupisan}  (and therefore \eqref{deterministic problem})
	and it is given by 
	\begin{equation}
	\label{semigroup solution u det 1}
	u (t,x)=   S_t g (x)+\int_0^t S_{t-s} f(s,x)\,ds, \quad t\in [0,T], \, x \in \mathbb{R}^d.
	\end{equation}
	To show the estimate \eqref{ee-perturbations}, we start from \eqref{semigroup solution u det 1}, use the bound \eqref{semigroupp estimate 2} together with \eqref{eModt} and  obtain
	\[\begin{split}
	\|u(t,\cdot)\|_{L^2(\mathbb R^d)} &\leq  M(t) \|g(\cdot)\|_{L^2(\mathbb R^d)} + \int_0^t M(t-s) \|f(s,\cdot)\|_{L^2(\mathbb R^d)}\,ds \\
	& \leq M(t) \left( \|g(\cdot)\|_{L^2(\mathbb R^d)} + \int_0^t \|f(s,\cdot)\|_{L^2(\mathbb R^d)}\,ds\right), 
	\end{split}\]
	for all $t\in [0,T]$ with $M(t)$ as in \eqref{eModt}, which completes the proof. \hfill $\square$

\vspace{0.3cm}
\subsubsection{Regularization of distributions and very weak solutions}\label{SubSecReg}
In the sequel we want to allow the potential $q$ to be a strongly irregular function or a distribution.
A common approach in the functional analysis to deal with irregularities is to regularize them.
The regularization can be always   achieved via convolutions with mollifiers and regularization families.

A smooth function $\varphi$ is a \emph{mollifier} if $\varphi \in \D (\mathbb{R}^d)$,  $\varphi \geq 0$ and $\int_{\mathbb{R}^d} \varphi (x)\,dx=1$. The \emph{mollifying   net}  $(\varphi_{\varepsilon})_{\varepsilon\in (0,1]}$  is 
\begin{equation} 
\label{molifajer}
\varphi_{\varepsilon} (x)= \dfrac{1}{(l(\varepsilon))^d} \, \varphi \left( \dfrac{x}{l(\varepsilon)}\right) \in \D(\mathbb{R}^d),
\end{equation}
where $l $ is a positive function and $l(\varepsilon)\to 0$, as $\varepsilon\to  0$. 

For a given distribution $q\in\D'(\R^d)$ a
regularization via convolution is the process in which we convolve  $q$ with the mollifying   net  $(\varphi_\varepsilon)_{\varepsilon \in (0,1]}$ and obtain the \emph{regularizing net} $(q_\varepsilon)_{\varepsilon\in (0,1]}$  of smooth functions  with compact supports 
\begin{equation}\label{reg_q}
q_\varepsilon(x):=q\ast \varphi_\varepsilon (x), \quad \varepsilon\in (0,1],
\end{equation}
which converges to $q$ in $\D'(\R^d)$ as $\varepsilon\to  0$. 

In order to work with regularizing nets in (S)PDEs, and to control the rate of their convergence we introduce  nets of moderate families in the following way.
\begin{de}[{\bf Moderate nets}] 
	\label{def: moderate nets} 
	Let $(\mathcal{B},\| \cdot\|_{\mathcal{B}})$ be a Banach space. A net of elements $(b_\varepsilon)_{\varepsilon \in (0,1]}$ in $\mathcal B$ is called $\mathcal{B}$-moderate if there exist $N \in \mathbb{N}_0$ and $C>0$ such that for each $\varepsilon \in (0,1]$ $$\|b_{\varepsilon}\|_{\mathcal{B}} \leq C \varepsilon^{-N}.$$	
\end{de}

There are other ways to define moderateness of which the most general one  is to define moderate nets  via families of seminorms as it was done in \cite{GR15}.   
Moderate nets are widely used in the theory of Colombeau algebras, along with the other notion of \emph{negligiblity} here also used for  defining uniqueness, see below Section \ref{SubsecConsUni}. We remark  that in the very weak solution concept it is  \emph{not} required derivatives to be moderate, as it is mostly the case when one applies the Colombeau theory of generalized functions in solving (S)PDEs  \cite{MO1992,MO2003}. For more details on Colombeu theory we refer to \cite{Colmbeau1984,GKOS2001,HO2009,HKO2013,NPS1998,MO1992}. 
Here using the very weak solution approach we choose nets more freely without caring whether the underlying space is an algebra or not, see for example \cite{GR15,MRT2019a,MRT2019b,RT17a,RT17b}.

In the following lemma we summarize some of the important results needed in sequel. For the details we refer  \cite{GKOS2001,MO1989,MO1992}.  
\begin{lem}
	\label{LemmaModerateness}
	Let $q\in \E'(\R^d)$.  
	\begin{itemize}
		\item[(a)] 
		The net $(q_\varepsilon)_{\varepsilon\in (0,1]}$ defined in \eqref{reg_q} is $L^\infty(\mathbb R^d)$- moderate. 
		\item[(b)]  
		There exists a mollifying   net  such that the net 
		$(q_\varepsilon)_{\varepsilon\in (0,1]}$ defined in \eqref{reg_q} is 
		$L^\infty(\mathbb R^d)$-log-type  moderate, i.e.,
		there exists $N>0$ so that it holds 
		\begin{equation*}
		\label{q_log-type} 
		\|q_{\varepsilon}\|_{L^{\infty}(\mathbb R^d)}  \leq N\, \log {\frac 1\varepsilon}. 
		\end{equation*}
	\end{itemize}
\end{lem}

{\bf Proof.} 
$(a)$ 
	By distribution structure theorem every compactly supported distribution is  of finite order and can be represented as a finite sum of derivatives of a continuous functions \cite{vladimirov1976generalized}. Therefore, i.e., there exist $k \in \mathbb{N}$ and $f \in C_0 (\R^d)$ such that $q= \sum_{|\alpha| \leq k} \partial^\alpha f.$ It holds
	\begin{align*}
	q\ast \varphi_\varepsilon (x) & = \sum_{|\alpha|\leq k} \partial^\alpha f \ast \varphi_\varepsilon (x) = \sum_{|\alpha|\leq k}  f  \ast \partial^\alpha \varphi_\varepsilon (x)\\
	&=\sum_{|\alpha|\leq k} \int f (y) \partial^{\alpha}_x \varphi_\varepsilon (x-y)\,dy\\
	&= \sum_{|\alpha|\leq k} \int_{\supp (\varphi_\varepsilon)} f (y) l(\varepsilon)^{-d-|\alpha|}\partial^{\alpha}_x \varphi \left( \frac{x-y}{l(\varepsilon)}\right)\,dy\\
	& \leq (l(\varepsilon))^{-k}   \sum_{|\alpha|\leq k}  \int f (x-yl(\varepsilon))  \partial^{\alpha}_y \varphi (y)\, dy.
	\end{align*} 
	Therefore, for each $\varepsilon \in (0,1]$ we have  $$\|q_{\varepsilon}\|_{L^\infty(\R^d)} = \|q\ast \varphi_\varepsilon \|_{L^\infty(\R^d)}\leq C (l(\varepsilon))^{-k}, 
	$$ with $C= \|f\|_{L^\infty(\R^d)} \sum_{|\alpha| \leq k} \int | \partial^{\alpha}_y \varphi (y)|\,dy.$
	If we choose $l(\varepsilon)=\varepsilon$ and $N=|\alpha|$, it follows that the net  $(q_\varepsilon)_{\varepsilon\in(0,1]}$ is $L^\infty(\R^d)$-moderate in the sense of Definition \ref{def: moderate nets}.
	
	$(b)$  
	By the rescaling mollification  process described in \cite[Prop. 1.5]{MO1989} it follows that every distribution of finite order can be regularized such that the obtained regularization net is of $L^{\infty}(\R^d)$-log-type. \hfill $\square$

\begin{rmk}\label{RemModerateness}(i) {Note that if $q \in L^{\infty} (\R^d)$ and the regularizing net $(q_\varepsilon)_{\varepsilon \in (0,1]}$ is defined in \eqref{reg_q}, then there exists $C>0$ such that for each $\varepsilon \in (0,1]$
		$$\|q_\varepsilon \|_{L^\infty (\R^d)} \leq C,$$ i.e., the regularizing net $(q_\varepsilon)_{\varepsilon \in (0,1]}$   is $L^\infty (\R^d)$-moderate in the sense of Definition \ref{def: moderate nets} with $N=0.$}
	
	(ii)  {Results of  Lemma \ref{LemmaModerateness} can be locally extended to  $q\in \D'(\R^d)$. This can be done by introducing an open covering of $\R^d$, corresponding cut-off functions and partition of unity subordinated to the covering. For the details we refer to \cite{GKOS2001}. }
\end{rmk}

In a similar way as it was done in \cite{ARST2020c} we define a very weak solution to the deterministic parabolic initial value  problem \eqref{deterministic problem} with potential {$q\in \D'(\R^d)$}  and then prove its existence, uniqueness and compatibility with the  solution obtained in  Theorem \ref{theorem estimates L infinity deterministic}.

\begin{de}[{\bf Very weak solutions for deterministic parabolic initial value  problem}] 
	\label{def: very weak solution}
	A net $(u_\varepsilon)_{\varepsilon\in (0,1]} \subset C([0,T]; D) \subseteq C([0,T]; L^2(\R^d))$ 
	is  \emph{very weak solution} to the problem \eqref{deterministic problem}
	if there exists a	regularizing net of smooth functions $(q_\varepsilon)_{\varepsilon\in (0,1]}$ of the potential {$q\in \D'(\R^d)$} such that  for every $\varepsilon\in (0,1]$, $u_\varepsilon$ is a solution to
	\begin{equation}
	\label{regularized deterministic problem}
	\begin{split}
	\left( \frac{\partial}{\partial t}  -  {\L} \right)  u(t, x) + {q_\varepsilon (x)} \cdot  u(t, x) & = f(t, x ), \quad t\in (0,T], \,   \, x\in  \mathbb R^d,   
	\\
	u(0, x) &= g(x),
	\end{split}
	\end{equation} 
	and $(u_\varepsilon)_{\varepsilon \in (0,1]}$ is $C([0,T]; L^2(\R^d))$-moderate.
\end{de}

We note that  \eqref{regularized deterministic problem}  gives the net of the problems which we will often referred as the \emph{net of  regularized problems} corresponding to \eqref{deterministic problem}.

\begin{thm}\label{thm extended Rozha}
	Let the operator $\L$ be as in (H1), the force term $f\in C([0,T]; L^2(\R^d)) $ a Lipschitz continuous function with respect to $t$, and the initial condition $g\in D$. 
	Let the potential {$q\in \D'(\R^d)$}. Then,  the  problem \eqref{deterministic problem} has a very weak solution in sense of Definition \ref{def: very weak solution}. 
\end{thm}

{\bf Proof.} 
	We start by regularizing  the potential  {$q\in \D'(\R^d)$}  by  means of convolution with a  mollifying net   $(\varphi_\varepsilon)_{\varepsilon\in (0,1]}$ given by \eqref{molifajer}, where $\varphi$ is chosen so that the obtained  net of smooth functions $q_{\varepsilon}(x) = q \ast\varphi_{\varepsilon}(x)$, $\varepsilon \in (0,1]$,  is of log-type (possible by Lemma \ref{LemmaModerateness} (b), and the Remark \ref{RemModerateness} (ii)),   i.e., there exists $N_q>0$ so that
	\begin{equation}
	\label{reg q varepsilon} 
	\|q_{\varepsilon}\|_{L^{\infty}(\mathbb R^d)}  \leq N_q \, \log {\frac 1\varepsilon}. 
	\end{equation}
	For each $\varepsilon \in(0,1]$, since $q_\varepsilon \in L^\infty (\R^d)$ by Theorem \ref{theorem estimates L infinity deterministic} there exists solution $u_\varepsilon \in C([0,T]; D)$ to the regularized problem  \eqref{regularized deterministic problem} 
	satisfying  the   estimate \eqref{ee-perturbations}, i.e., 
	$$\|u_\varepsilon (t,\cdot)\|_{L^2(\mathbb R^d)} \leq   Me^{(w+M \|q_\varepsilon (\cdot)\|_{L^\infty (\R^d)}) t} \left( \|g(\cdot)\|_{L^{2}(\mathbb R^d)} +\int_0^t \|f(s,\cdot)\|_{L^2(\mathbb R^d)} \,ds\right) 
	$$for $t\in [0,T]$.  We obtain the net 
	$(u_\varepsilon)_{\varepsilon\in(0,1]} \subset C([0,T]; D) \subseteq  C([0,T]; L^2(\R^d)),$
	which will represent a very weak solution in the sense of Definition  \ref{def: very weak solution} only if it is $C([0,T]; L^2(\R^d))$-moderate. 
	By  \eqref{reg q varepsilon} we further obtain 
	\begin{equation*}
	\begin{split}
	\|u_\varepsilon (t,\cdot)\|_{L^2(\mathbb R^d)} & \leq  Me^{(w+M \|q_\varepsilon (\cdot)\|_{L^{\infty}(\mathbb R^d)}) t} \left( \|g(\cdot)\|_{L^{2}(\mathbb R^d)} +\int_0^t \|f(s,\cdot)\|_{L^2(\mathbb R^d)} \,ds\right) \\
	& \leq Me^{wt} e^{- M \, N_q\,t \, \log \varepsilon}  \left( \|g(\cdot)\|_{L^{2}(\mathbb R^d)} +\int_0^t \|f(s,\cdot)\|_{L^2(\mathbb R^d)} \,ds\right)\\
	&=  Me^{wt}  \left( \|g(\cdot)\|_{L^{2}(\mathbb R^d)} +\int_0^t \|f(s,\cdot)\|_{L^2(\mathbb R^d)} \,ds\right) \,\,  \varepsilon^{-N}\\
	& = c(t) \,   \varepsilon^{-N}, 
	\end{split}
	\end{equation*}
	where $c(t)=Me^{wt}  ( \|g(\cdot)\|_{L^{2}(\mathbb R^d)} +\int_0^t \|f(s,\cdot)\|_{L^2(\mathbb R^d)} \,ds)$ and $N= M \,  N_q\,t$ are  positive constants for each $t\in (0,T]$. 
	Finally, we have
	\[\|u_\varepsilon\|_{C([0,T]; L^2(\R^d))} =  
	\sup_{t\in[0,T] } ( \|u_{\varepsilon}(t, \cdot) \|_{L^2(\R^d)} )
	\leq \sup_{t\in[0,T] } ( c(t)\, \varepsilon^{-N}) \leq C(T) \,  \varepsilon^{-N}.\]
	Therefore, there exist $N\in \N_0$ and $C=C(T) > 0$  such that 
	\begin{equation*}
	\|u_\varepsilon\|_{C([0,T]; L^2(\R^d))} \leq C \, \varepsilon^{-N}
	\end{equation*} holds, i.e., $(u_\varepsilon)_{\varepsilon\in(0,1]} $   is a
	$C([0,T]; L^2(\R^d))$-moderate net. \hfill $\square$

\begin{rmk}
	\label{RemLogType}  
	\begin{itemize}
		\item[(i)]	Note that in order to have a moderate solution family, we have to choose a regularization of $q$ to be   of log-type.   
		\item[(ii)]  If the net $(u_\varepsilon)_{\varepsilon \in (0,1]}$ is $C([0,T]; L^2(\R^d))$-moderate, then the net \linebreak $(u_\varepsilon(t,\cdot) )_{\varepsilon \in (0,1]}$ is $L^2(\R^d)$-moderate. Indeed, there exist $N \in \mathbb{N}_0$ and $C>0$ such that
		$$\|u(t,\cdot)\|_{L^2(\R^d)} \leq \sup_{t \in [0,T]} \|u(t,\cdot)\|_{L^2(\R^d)} = \|u_\varepsilon\|_{C([0,T]; L^2(\R^d))} \leq C \varepsilon^{-N}.$$
		Along the proof of the previous theorem, one can see that the oposite is also true. 
	\end{itemize}
\end{rmk} 
Theorem \ref{thm extended Rozha}  extends the results from \cite{ARST2020c}, where the authors considered the heat equation with singular potential, but without a force term. Namely, they studied a special case of \eqref{deterministic problem} with $ \L=\Delta$ and $f=0$ and showed existence and uniqueness of a very weak solution as well as consistency of  the obtained solution with the classical solution. In what follows we will obtain the same for our problem \eqref{deterministic problem}. The semilinear heat equation with singular potential and a non-zero force term, i.e., the problem \eqref{deterministic problem} with $ \L=\Delta$, is considered in \cite{NPR2005} using semigroup approach in the framework of Colombeau theory. %

\vspace{0.2cm}
\subsubsection{Questions on uniqueness and consistency with classical weak solutions}\label{SubsecConsUni}

The question of uniqueness for the very weak solution  obtained in Theorem \ref{thm extended Rozha} can be treated in different manners. 
If we assume that $(u_\varepsilon)_{\varepsilon\in(0,1]}$ and $(\tilde{u}_\varepsilon)_{\varepsilon\in(0,1]}$ are two very weak solutions obtained in Theorem \ref{thm extended Rozha}, according to Definition \ref{def: very weak solution} it means that there exist regularizing nets $(q_\varepsilon)_{\varepsilon\in(0,1]}$ and $(\tilde{q}_\varepsilon)_{\varepsilon\in(0,1]}$ such that for each $\varepsilon\in(0,1]$, $u_\varepsilon$ and $\tilde{u}_\varepsilon$ satisfy the net of regularized problems \eqref{regularized deterministic problem} with $q_\varepsilon$ and $\tilde{q}_\varepsilon$ as  potentials respectively, and both nets  $(u_\varepsilon)_{\varepsilon\in(0,1]}$ and $(\tilde{u}_\varepsilon)_{\varepsilon\in(0,1]}$ are $C([0,T]; L^2(\R^d))$-moderate. If it happens that $q_\varepsilon = \tilde{q}_\varepsilon$, for all $\varepsilon\in(0,1]$ (i.e., these are two same regularizing nets), then it is straightforward to show that it implies $u_\varepsilon = \tilde{u}_\varepsilon$,  for all $\varepsilon\in(0,1]$.
Indeed, if we denote $v_\varepsilon := u_\varepsilon - \tilde{u}_\varepsilon$, $\varepsilon \in (0,1],$ then,  $(v_\varepsilon)_{\varepsilon\in(0,1]}$ is a very weak solution to  the  homogeneous problem corresponding to \eqref{deterministic problem}, i.e. for each $\varepsilon\in(0,1]$ it holds
\begin{equation*}
\begin{split}
\left(\frac{\partial}{\partial t}  - \mathcal L \right) v_\varepsilon(t, x) + {q_\varepsilon(x)} v_\varepsilon (t, x) &= 0, 
\\
v_\varepsilon(0, x) &=0 . 
\end{split}
\end{equation*}
Applying Theorem \ref{theorem estimates L infinity deterministic} with $f=g=0$, and using the estimate \eqref{ee-perturbations} it follows that for each $\varepsilon \in (0,1],$ $v_{\varepsilon} = 0$,   and therefore $u_{\varepsilon} = \tilde{u}_{\varepsilon}$, $\varepsilon \in (0,1]$.

In the case when regularizations $(q_\varepsilon)_{\varepsilon\in(0,1]}$ and $(\tilde{q}_\varepsilon)_{\varepsilon\in(0,1]}$ are not the same, then uniqueness holds in the following sense.
\begin{thm}[{\bf Uniqueness}] \label{thm:uniqueness deterministic}
	Let $(q_\varepsilon)_{\varepsilon\in(0,1]}$ and $(\tilde{q}_\varepsilon)_{\varepsilon\in(0,1]}$ be two regularizing nets of {$q\in\D'(\R^d)$} so that for all $n\in \N$ there exists $c>0$ such that
	\begin{equation}
	\label{q-IS-negligible}
	\|q_\varepsilon-\tilde{q}_\varepsilon\|_{L^\infty(\R^d)}\leq c \,\varepsilon^n.
	\end{equation}
	Let  $(u_\varepsilon)_{\varepsilon\in(0,1]}$ and $(\tilde{u}_\varepsilon)_{\varepsilon\in(0,1]}$ be two very weak solutions of the problem \eqref{deterministic problem} with respective regularizing nets $(q_\varepsilon)_{\varepsilon\in(0,1]}$ and $(\tilde{q}_\varepsilon)_{\varepsilon\in(0,1]}$ of the potential {$q\in \D'(\R^d)$}. Then,
	for all $k\in \N$ there exists $C>0$ such that
	\begin{equation*}
	\|u_\varepsilon-\tilde{u}_\varepsilon\|_{_{C([0,T]; L^2(\R^d))}}\leq C \, \varepsilon^k.
	\end{equation*}  
\end{thm}

{\bf Proof.} 
	By assumptions, for each $\varepsilon \in (0,1]$, $u_\varepsilon$ and $\tilde{u}_\varepsilon$ satisfy the regularized problem \eqref{regularized deterministic problem} with $q_\varepsilon$ and $\tilde{q}_\varepsilon$ as  potentials respectively, and both nets  $(u_\varepsilon)_{\varepsilon\in(0,1]}$ and $(\tilde{u}_\varepsilon)_{\varepsilon\in(0,1]}$ are $C([0,T]; L^2(\R^d))$-moderate. In the view of Remark \ref{RemLogType} (i) both nets $(q_\varepsilon)_{\varepsilon\in(0,1]}$ and $(\tilde{q}_\varepsilon)_{\varepsilon\in(0,1]}$ are of log-type. 
	Then it holds 
	\begin{equation*}
	\begin{split}
	\frac{\partial}{\partial t} (u_{\varepsilon} - \tilde{u}_{\varepsilon})(t, x) - \mathcal{L}  (u_{\varepsilon} - \tilde{u}_{\varepsilon})(t, x) + {q_\varepsilon (x)} \cdot  (u_{\varepsilon} - \tilde{u}_{\varepsilon})(t, x) & = f_{\varepsilon}(t, x ),    
	\\
	(u_{\gamma,\varepsilon} - \tilde{u}_{\gamma,\varepsilon})(0, x) =0, 
	\end{split}
	\end{equation*}
	with $$f_{\varepsilon} (t,x) = (\tilde{q}_{\varepsilon} - q_\varepsilon)(x) \tilde{u}_{\varepsilon} (t,x).$$ 
	Since $(\tilde{u}_\varepsilon)_{\varepsilon\in (0,1]}$ is a very weak solution,  it is $C([0,T]; L^2(\R^d))$-moderate, but also $(\tilde{u}_{\varepsilon}(t,\cdot))_{\varepsilon\in (0,1]}  $ is $L^2(\R^d)$--moderate (see Remark \ref{RemLogType} (ii)), i.e.,  there exist $N \in \mathbb{N}_0$ and $C_1>0$ such that 
	\begin{equation} 
	\label{koeficijenti moderate}
	\|\tilde{u}_{\varepsilon} (t,\cdot) \|_{L^2(\R^d)} \leq C_1 \, \varepsilon^{-N}.
	\end{equation}
	Using the estimate \eqref{ee-perturbations} together with \eqref{q-IS-negligible}, and \eqref{koeficijenti moderate},  and the fact that the regularizing net $(q_\varepsilon)_{\varepsilon\in (0,1]}$ is of log-type, we obtain 
	\begin{equation*}
	\begin{split}
	\|(u_{\varepsilon} &- \tilde{u}_{\varepsilon})(t,\cdot)\|_{L^2(\mathbb R^d)}  \leq  Me^{wt}e^{M\|q_{\varepsilon} (\cdot)\|_{L^\infty(\R^d)} t} \int_{0}^t \|(\tilde{q}_\varepsilon - q_{\varepsilon})(\cdot)\tilde{u}_{\varepsilon}(s, \cdot)\|_{L^2(\mathbb R^d)}\,ds \\
	& \leq  Me^{wt}e^{M\|q_{\varepsilon} (\cdot)\|_{L^\infty(\R^d)} t} \|(\tilde{q}_\varepsilon - q_{\varepsilon})(\cdot)\|_{L^{\infty}(\mathbb{R}^d)} \int_0^t \|\tilde{u}_{\varepsilon}(s, \cdot)\|_{L^2(\mathbb R^d)}\,ds \\
	&\leq  Me^{wt}\varepsilon^{-M N_q t} \, C \varepsilon^n \, \,  \int_0^t C_1 \varepsilon^{-N}\,ds  \\
	& =  Mt \, e^{wt} C_1\, \varepsilon^{n-N-M N_q t}  \\
	&=  c(t) \, \varepsilon^{n-N-M N_q t} 
	\end{split}
	\end{equation*} 
	for arbitrary $n.$  We choose $n> N_0:=N+M N_q T$ to obtain 
	$$	\|u_{\varepsilon}-\tilde{u}_{\varepsilon}\|_{C([0,T]; L^2(\R^d))} =
	\sup_{t\in[0,T] }  \|(u_{\varepsilon}-\tilde{u}_{\varepsilon})(t, \cdot) \|_{L^2(\R^d)} \leq  C(T) \,  \varepsilon^{k}$$	
	for all $k\in \mathbb N$.   \hfill $\square$

The question whether the notion of very weak solution obtained in Theorem \ref{thm extended Rozha} for the case when the potential is bounded regular function coincides  with the  solution guaranteed by Theorem \ref{theorem estimates L infinity deterministic} is treated in the following theorem.  

\begin{thm}[{\bf Consistency}] \label{thm:consistency deterministic}
	Let the operator $\L$ be as in (H1), the force term $f\in C([0,T]; L^2(\R^d))$ is a Lipschitz continuous function with respect to $t$, the initial condition $g\in D$,  and let the potential $q\in L^\infty(\mathbb R^d)\cap C(\mathbb R^d)$.
	
	Let $v\in C([0,T];D)$ be the weak solution  to the parabolic initial value problem \eqref{deterministic problem} obtained by Theorem \ref{theorem estimates L infinity deterministic} and let $(u_\varepsilon)_{\varepsilon\in (0,1]} \subset C([0,T]; L^2(\R^d))$  be the very weak solution to  \eqref{deterministic problem} obtained by Theorem \ref{thm extended Rozha}. Then, 
	\begin{equation*}
	\label{Converg}
	\|u_\varepsilon - v \|_{C([0,T]; L^2(\R^d))} \to 0 \qquad \text{as} \qquad \varepsilon \to 0.
	\end{equation*}
\end{thm} 

{\bf Proof.} 
	Let $(q_\varepsilon)_{\varepsilon \in (0,1]}$ be a regularizing net \eqref{reg_q} of given $q\in L^\infty(\mathbb R^d)\cap C(\mathbb R^d)$. 
	Then $(q_\varepsilon)_{\varepsilon \in (0,1]}$ 
	converges to $q$ uniformly, i.e., $ \|{q} - q_\varepsilon\|_{L^\infty(\R^d)}\to 0$, as $\varepsilon\to 0$.
	By  the assumptions, for every $\varepsilon\in (0,1]$, $u_\varepsilon$ satisfies \eqref{regularized deterministic problem} with $q_\varepsilon$ from the above regularization, while $v$ satisfies \eqref{deterministic problem}.  
	We subtract \eqref{regularized deterministic problem} and \eqref{deterministic problem} and obtain
	\begin{align*}
	\frac{\partial}{\partial t} (u_{\varepsilon}  - v)  -  \mathcal{L}   (u_{\varepsilon}  -v) + {q_\varepsilon} \cdot u_{\varepsilon}  - q \cdot   v & = 0, \\
	(u_{\varepsilon}  -v)|_{t=0}&= 0. 
	\end{align*}
	After adding and subtracting $q_\varepsilon v$, and regrouping the factors,  the previous equation transforms to 
	\[\begin{split}
	\left(\frac{\partial}{\partial t} -\mathcal L + q_\varepsilon \cdot Id\right) \,  (u_{\varepsilon}  - v)    & =  (q- q_\varepsilon) \, v,  \\
	(u_{ \varepsilon}  -v)|_{t=0} &= 0. 
	\end{split}\]
	Hence, for all $\varepsilon$ the difference $u_{\varepsilon}  -v$   satisfies \eqref{deterministic problem} with zero initial condition,  in which $q$ is replaced by $q_\varepsilon$ and $f$ is replaced by $(q- q_\varepsilon) \, v$. Thus, by Theorem \ref{theorem estimates L infinity deterministic},  the estimate \eqref{ee-perturbations} and Remark \ref{RemModerateness} (i) the difference $u_{\varepsilon}  -v$ satisfies
	\begin{equation} 
	\begin{split}
	\|(u_{\varepsilon}  -v)(t,\cdot)\|_{L^2(\mathbb R^d)} &\leq   M e^{wt}e^{M \|q_\varepsilon (\cdot)\|_{L^\infty (\R^d)}t} \int_0^t \|(q-q_\varepsilon)(\cdot) v(s,\cdot)\|_{L^2(\mathbb R^d)}\,ds \nonumber \\ 
	&\leq   M e^{wt}e^{M C t}  \|(q-q_\varepsilon)(\cdot) \|_{L^\infty (\mathbb{R}^d)} \int_{0}^t \|v(s,\cdot)\|_{L^2(\mathbb R^d)} \,ds \nonumber \\  
	&\leq \tilde{C}(t)  \|(q-q_\varepsilon)(\cdot) \|_{L^\infty (\mathbb{R}^d)}, 
	\end{split}
	\end{equation} 
	where  $\tilde{C}(t):= M e^{(w+MC)t} \int_{0}^t \|v(s,\cdot)\|_{L^2(\mathbb R^d)} \,ds $  is a finite  constant. Finally, since $ \|{q} - q_\varepsilon\|_{L^\infty(\R^d)}\to 0$, as $\varepsilon\to 0$ we have also that
	$\|u_{\varepsilon}  - v\|_{C([0,T];L^2(\R^d))} = \sup_{t \in[0,T]} \| (u_{\varepsilon}  - v)(t,\cdot)\|_{L^2(\R^d)} $
	tends to $0$ as $\varepsilon \to 0$. \hfill $\square$

\vspace{0.3cm}
\section{A very weak solutions for the stochastic parabolic equations with singular space depending potentials} \label{sec4}

We are now ready to analyse the stochastic parabolic problem  \eqref{Eq: stochastic evolution} for {$q\in\D'(\R^d)$} combining   chaos expansion representations and the concept of very weak solutions. 	In the first step, we proceed as in Theorem \ref{stochastic weak sol} and represent all stochastic processes appearing in  \eqref{Eq: stochastic evolution} in their chaos expansion forms. We also represent a solution in the chaos expansion form \eqref{ChaosExpansion}, i.e.,  
\begin{equation}
\label{sol U}
U(t, x, \omega) = \sum\limits_{\gamma\in \mathcal I} u_\gamma(t, x) \, H_\gamma(\omega), \qquad t\in[0,T], \,\, x\in \mathbb R^d, \,\, \omega\in \Omega.
\end{equation}
After  equalizing corresponding coefficients in the Fourier-Hermite polynomial basis, as in Theorem  \ref{stochastic weak sol},  we obtain deterministic problems, but in this case now with  irregular potential {$q \in \mathcal{D}'(\R^d)$}. For each $\gamma \in \mathcal{I}$ the existence of a very weak solution to the obtained deterministic  PDE  is guaranteed by Theorem \ref{thm extended Rozha} which determines the coefficients $u_\gamma$, $\gamma \in \mathcal{I}$ of the solution \eqref{sol U}. Since very weak solution is a net, with this procedure we  obtain a net  $(U_\varepsilon)_{\varepsilon \in (0,1]}$, for which we have to prove that it belongs to a Kondratiev-type  space of generalized stochastic processes. 
The notion for moderate nets of generalized stochastic processes appears naturally as follows.

\begin{de}[{\bf Moderate nets of generalized stochastic processes}] 
	\label{def: moderate nets of processes} 
	Let $(\mathcal{B},\| \cdot\|_{\mathcal{B}})$ be a Banach space. 	Let  $(U_\varepsilon)_{\varepsilon \in (0,1]}\subset \mathcal{B} \otimes (S)_{-1}$ be a  net of generalized stochastic processes $U_\varepsilon \in \mathcal{B} \otimes (S)_{ -p_\varepsilon}$,  $p_\varepsilon \geq 0$,  and for every $\varepsilon \in (0,1]$ 
	\begin{equation} 
	\label{U epsilon net chaos exp}
	U_\varepsilon = \sum_{\gamma\in \mathcal I} (u_\gamma)_{\varepsilon} \, H_\gamma.
	\end{equation} 
	Assume $p := \sup\limits_{\varepsilon \in (0,1]} \, p_\varepsilon$ exists. If for all $\gamma \in \mathcal I$ the  coefficients $[(u_\gamma)_{\varepsilon}]_{\varepsilon \in (0,1]}$  are $\mathcal{B}$-moderate nets, then $(U_\varepsilon)_{\varepsilon \in (0,1]}$ is  called $\mathcal{B} \otimes (S)_{-p}$-\emph{moderate}.  
\end{de}
\begin{rmk}
	\label{Remark2}
	\begin{itemize}
		\item[(i)] If the supremum $p = \sup\limits_{\varepsilon \in (0,1]} \, p_\varepsilon$ does not exist, one can define moderateness in the Frech\'et space  $\mathcal{B} \otimes (S)_{-1}$ with seminorms estimates as it was done in \cite{Garreto:05}. For simplicity and the presentation of the method, here we work only with Banach spaces. 
		\item[(ii)] If $(U_\varepsilon)_{\varepsilon \in (0,1]}$ is $\mathcal{B} \otimes (S)_{ -p}$-moderate then it is also  $\mathcal{B} \otimes (S)_{ -r}$-moderate for each $r\geq p$.
		\item[(iii)] Note that if $p=0$ then the net $(U_\varepsilon)_{\varepsilon\in (0,1]}$ is  $\mathcal{B}\otimes L^2 (\mu)$-moderate.
	\end{itemize}
\end{rmk}

One can proceed differently, starting with the regularization of distribution $q$,  one can consider a net of regularized stochastic parabolic equations, then solve each of them by  Theorem \ref{stochastic weak sol}, which will result in the net of generalized stochastic processes.  
To have moderateness for  this net  a different notion for moderate nets appears as natural.  For a better insight,  we state and prove the following lemma, which assures that two approaches are equivalent.

\begin{lem}\label{Lemma solution concepts are equivalent} 
	Assume    $(U_\varepsilon)_{\varepsilon\in (0,1]}\subset \mathcal{B} \otimes (S)_{-1}$ is a net of generalized stochastic process $U_\varepsilon\in \mathcal{B} \otimes (S)_{ -p_\varepsilon}$, $p_\varepsilon\geq 0$. Assume that there exists $ p=\sup_{\varepsilon \in (0,1]} p_\varepsilon$, and $p > 1$.
	Then the net  $(U_\varepsilon)_{\varepsilon\in (0,1]}$  is  
	$\mathcal{B} \otimes (S)_{ -p}$-moderate 
	if and only if	
	there exist $N \in \mathbb{N}_0$  and $C>0$  such that 
	\begin{equation}
	\label{UepsMod}
	\|U_{\varepsilon}\|_{\mathcal{B}\otimes (S)_{ -p}}  \leq C \varepsilon^{-N}.
	\end{equation}  
\end{lem}

{\bf Proof.}
	If  $(U_\varepsilon)_{\varepsilon \in (0,1]}$ is $\mathcal{B}  \otimes (S)_{-p}$-moderate then  for all $\gamma \in \mathcal I$  the  coefficients $[(u_\gamma)_{\varepsilon}]_{\varepsilon \in (0,1]}$ in the chaos expansion \eqref{U epsilon net chaos exp} are $\mathcal{B} $-moderate,  
	i.e.,  there exist $N_1\in \mathbb N$ and $C_1>0$ such that 
	\begin{equation} 
	\label{part 2}
	\|(u_\gamma)_{\varepsilon}\|_{\mathcal{B} } \leq C_1 \, \varepsilon^{-N_1}. 
	\end{equation}
	Also, for each $\varepsilon \in (0,1]$,  $U_\varepsilon$   is a stochastic process in $\mathcal{B} \otimes (S)_{-p_\varepsilon }\subseteq \mathcal{B} \otimes (S)_{ -p}$, i.e., it holds
	\begin{equation} 
	\label{part 1}
	\|U_\varepsilon\|^2_{\mathcal{B}  \otimes (S)_{ -p}} = \sum_{\gamma\in \mathcal I} \|(u_\gamma)_{\varepsilon}\|^2_{\mathcal{B} }  \, (2\mathbb N)^{-p\gamma} < \infty.
	\end{equation}
	Combining  \eqref{part 2} and \eqref{part 1}  we obtain 
	\begin{equation*}
	\begin{split}
	\|U_\varepsilon\|^2_{\mathcal{B}  \otimes (S)_{ -p}} = 
	\sum_{\gamma\in \mathcal I} \|(u_\gamma)_{\varepsilon}\|^2_{\mathcal{B} }  \, (2\mathbb N)^{-p\gamma}
	\leq 
	C_1^2 \,  \varepsilon^{-2N_1} \,  \sum_{\gamma\in \mathcal I}  \, (2\mathbb N)^{-p\gamma} =  C_1^2  C_2 \, \varepsilon^{-2N_1},
	\end{split}
	\end{equation*}
	with $C_2=\sum_{\gamma\in \mathcal I}  \, (2\mathbb N)^{-p\gamma} $ being finite for $p>1$ by \eqref{2N konv p>1}.
	If we take $C=C_1^2C_2$ and $N=2N_1$ we finish the first part of the proof. 
	
	Conversely, suppose that there exist $N \in \mathbb{N}_0$ and $C>0$ such that \eqref{UepsMod} holds with $  p = \sup_{\varepsilon \in (0,1]} \, p_\varepsilon$. 
	Suppose that the net $(U_\varepsilon)_{\varepsilon \in (0,1]}$ is not $\mathcal{B}  \otimes (S)_{ -p}$-moderate. This means that  for some $\gamma_0 \in \mathcal I$ the coefficients $(u_{\gamma_0})_{\varepsilon}$ are not $\mathcal{B}-$moderate, i.e.,  for all $N \in \mathbb{N}_0$ and $C>0$ it holds
	\begin{equation} 
	\label{non mod}
	\|(u_{\gamma_0})_{\varepsilon}\|_{\mathcal{B} }>C \varepsilon^{-N}. 
	\end{equation} 
	Therefore, using  \eqref{part 1} and \eqref{non mod} we obtain that for all $N \in \mathbb{N}_0$ and $C>0$ it holds
	\begin{eqnarray*}
		\|U_{\varepsilon}\|_{\mathcal{B} \otimes (S)_{-p}} & =& \sum_{\gamma \in \mathcal{I}} \|(u_\gamma)_{\varepsilon}\|_{\mathcal{B} }^2 (2\N)^{-p\gamma}\\
		& =& \sum_{\gamma \in \mathcal{I} \setminus \{\gamma_0\}} \|(u_\gamma)_{\varepsilon}\|_{\mathcal{B} }^2 (2\N)^{-p\gamma} + \|(u_{\gamma_0})_{\varepsilon}\|_{\mathcal{B} }^2 (2\N)^{-p\gamma_0}\\
		& > & \sum_{\gamma \in \mathcal{I} \setminus \{\gamma_0\}} \|(u_\gamma)_{\varepsilon}\|_{\mathcal{B} }^2 (2\N)^{-p\gamma} + C^2 \varepsilon^{-2N}(2\N)^{-p\gamma_0}\\
		& > & \tilde{C} \varepsilon^{-N},
	\end{eqnarray*}
	which leads to contradiction. \hfill  $\square$

Since $\mathcal{B} \otimes (S)_{-p_\varepsilon}$, $p_\varepsilon \geq 0$ are Banach spaces, according to Definition \ref{def: moderate nets} and Lemma \ref{Lemma solution concepts are equivalent}, moderateness  of generalized stochastic processes can be defined in the following way. 
\begin{de} \label{moderate gsp}
	Let $(\mathcal{B},\| \cdot\|_\mathcal{B})$ be a Banach space. 	Let  $(U_\varepsilon)_{\varepsilon \in (0,1]}\subset \mathcal{B} \otimes (S)_{-1}$ be a  net of generalized stochastic processes $U_\varepsilon \in \mathcal{B} \otimes (S)_{ -p_\varepsilon}$,  $p_\varepsilon \geq 0$, and let $  p = \sup\limits_{\varepsilon \in (0,1]} \, p_\varepsilon$ exists.
	The net of generalized stochastic processes 
	$(U_\varepsilon)_{\varepsilon \in (0,1]}$ is $\mathcal{B} \otimes (S)_{ -p}$-\emph{moderate}
	if	there exist $N \in \mathbb{N}_0$   and $C>0$  such that $$\|U_{\varepsilon}\|_{\mathcal{B}\otimes (S)_{ -p}}  \leq C \varepsilon^{-N}.$$  
\end{de}

\vspace{0.05cm}
\subsection{A stochastic very weak solution: Definition and existence} \label{subsec 3.1}
In this section we define  a stochastic very weak solution to the stochastic parabolic equation \eqref{Eq: stochastic evolution} with the potential {$q\in \D'(\R^d)$}, and prove its existence. 
In the following we use  $X:=C([0,T]; L^2(\R^d)) $ and recall that ${\cal X}(L^2(\R^d)):= X \otimes (S)_{-1}$.

\begin{de}[{\bf Solution concept 1}] \label{option 1} Let the potential $q$ be in  {$\mathcal D'(\mathbb R^d)$}.	A net  of stochastic processes $(U_{\varepsilon})_{\varepsilon\in (0,1]}$ in  $X \otimes (S)_{-1}$, which is $X \otimes (S)_{-p}$-moderate  in the sence of Definition \ref{def: moderate nets of processes}
	is  a {\rm very weak solution} to  the stochastic  parabolic  problem \eqref{Eq: stochastic evolution}	if  for each $\gamma\in \mathcal I$ the net  $[(u_{\gamma})_{\varepsilon}]_{\varepsilon\in (0,1]}$ of the coefficients in expansion \eqref{U epsilon net chaos exp}  is a very weak solution to 
	\begin{equation*}
	\begin{split}
	\left(\frac{\partial}{\partial t}  - \L \right)  \, u_\gamma(t, x) + {q(x)} \cdot u_\gamma(t, x) &= f_\gamma(x, t) \quad t\in (0,T], \, x\in \mathbb R^d,
	\\
	u_\gamma(0, x) &= g_\gamma(x).
	\end{split}
	\end{equation*} 
\end{de}

\begin{de}[{\bf Solution concept  2}]\label{option 2} Let the potential $q$ be in  {$\mathcal D'(\mathbb R^d)$}.	A net  of stochastic processes $(U_{\varepsilon})_{\varepsilon\in (0,1]}$ in  $X \otimes (S)_{-1}$
	is  a {\rm very weak solution} to  the stochastic  problem 
	\eqref{Eq: stochastic evolution}
	if  there exists an $L^\infty(\mathbb R^d)$-moderate  regularizing net $(q_\varepsilon)_{\varepsilon\in (0,1]}$ of the distribution $q$, such that for every $\varepsilon\in (0,1]$  the process $U_\varepsilon$ is a solution to the regularized stochastic parabolic  equation
	\begin{equation*}
	\label{SP reg}
	\begin{split}\small 
	\left( \frac{\partial}{\partial t}  - \mathcal{L} \right) \, U(t, x,\omega) + {q_\varepsilon (x)} \cdot  U(t, x,\omega) &= F(t, x,\omega),	\\
	U(0, x,\omega) &= G(x,\omega),
	\end{split}
	\end{equation*}
	and it is a moderate net of stochastic processes {in the sence of Definition \ref{moderate gsp}. }
\end{de}

\begin{rmk}
	\label{Remark3}
	We note that we work with double nets $(\gamma, \varepsilon)$, where $\gamma\in \mathcal I$ and $\varepsilon \in (0,1]$, i.e., for each $\gamma\in \mathcal I$ we have a sequence $(u_\gamma)_\varepsilon$, $\varepsilon \in (0,1]$. 
	In the solution concept 1  for each $\varepsilon\in (0,1]$ the stochastic process $U_\varepsilon$  is given by  
	\[U_\varepsilon = \sum_{\gamma\in \mathcal I} (u_\gamma)_\varepsilon  \, H_\gamma,\]
	while in the the solution concept 2	 for each $\varepsilon\in (0,1]$ the stochastic process $U_\varepsilon$ is given by 
	\[U_\varepsilon = \sum_{\gamma\in \mathcal I} (u_\varepsilon)_\gamma  \, H_\gamma.\]
	From Lemma \ref{Lemma solution concepts are equivalent} it follows that solution concept 1 and solution concept 2 are equivalent for $p>1$.  As the restriction $p\geq 0$ to   $p>1$ is not essential in the chaos expansion setting, in the following we will assume $p>1$ and use the equivalence of these two concepts. 
	In the sequel, we  will also write $u_{\gamma,\varepsilon}$ instead of $(u_\gamma)_\varepsilon$ or $(u_\varepsilon)_\gamma $.
\end{rmk}

\begin{thm}[{\bf Existence of a stochastic very weak solution}] \label{Thm 1: exist}
	Let the operator $\mathcal L$, the force term $F$ and the initial condition $G$ be such that the assumptions $(H1)$-$(H3)$ hold and let the  potential {$q\in \mathcal D'(\mathbb R^d)$}. Then,  the stochastic parabolic problem \eqref{Eq: stochastic evolution}	 
	has a  very weak solution  $(U_\varepsilon)_{\varepsilon\in (0,1]}$  in  $\mathcal{X}(L^2 (\R^d))$.
\end{thm}

{\bf Proof.} 
	Firstly, we write processes $F$ and $G$  in their chaos expansion forms as in  Theorem \ref{stochastic weak sol}, i.e.,  $F$ is of the form \eqref{F chaos}  such that  the  condition \eqref{force term uslov konv} holds for some  $p_1\geq 0$,  and $G$ is of the form \eqref{G chaos} such that the condition \eqref{poc uslov konv} holds for some $p_2\geq 0$.  
	We also assume that  the unknown process is given in the form  \eqref{sol U}. 
	After applying the chaos expansion method, for every $\gamma \in \mathcal I$ we obtain the equation	
	\begin{equation}
	\label{problemm Pgamma}
	\begin{split}
	\left(\frac{\partial}{\partial t} -  \L  \right)u_\gamma(t, x) + q(x) \cdot u_\gamma(t, x) &= f_\gamma(x, t), \\
	u_\gamma(0, x) &= g_\gamma(x),
	\end{split}
	\end{equation}
	where {$q\in \mathcal D'(\mathbb R^d)$}.
	From the assumptions $(H1)$-$(H3)$  we have for all $\gamma\in \mathcal I$ that  $f_\gamma \in C([0,T];L^2(\R^d))$ are  Lipschitz continuous functions with respect to $t$ and $g_\gamma \in D$.  By Theorem \ref{thm extended Rozha}, for each $\gamma\in \mathcal I$ there exists a very weak solution to \eqref{problemm Pgamma}  i.e., a $X$-moderate net  
	$ [(u_\gamma)_{\varepsilon}]_{\varepsilon \in (0,1]} = (u_{\gamma, \varepsilon})_{\varepsilon \in (0,1]}$  in $C([0,T],D )$, satisfying \eqref{problemm Pgamma} and 
	there exists $N>0$ such that  
	\begin{equation}
	\label{estimate u gama epsilon nets}
	\|u_{\gamma, \varepsilon} (t,\cdot)\|_{L^2(\R^d)} \leq Me^{wt}  \left( \|g_\gamma (\cdot)\|_{L^{2}(\R^d)} +\int_0^t \|f_\gamma (s,\cdot)\|_{L^2(\R^d)} \,ds\right)  \varepsilon^{-N}.
	\end{equation}
	Note that for each $\gamma\in \mathcal I$, the  chosen constant $N= N(M, T, w, N_q)$  does not depend on $\gamma$. 
	Finally, we are going to show that 
	for every $\varepsilon \in (0,1]$ a process defined by  
	\begin{equation*}
	\label{Stoch net procesi}
	U_{\varepsilon} (t,x,\omega): = \sum_{\gamma \in \mathcal{I}} u_{\gamma, \varepsilon} (t,x) \, H_{\gamma}(\omega),
	\end{equation*}	
	is a  generalized stochastic process in $\mathcal{X}(L^{2}(\R^d))$ which is moderate,  i.e.,  we have to show that for each $\varepsilon \in (0,1]$ there exists $p_\varepsilon>0$ such that  $U_\varepsilon \in  X \otimes (S)_{ -p_\varepsilon}$ and that \eqref{UepsMod} holds.  
	We are going  to  prove that for all $\varepsilon\in (0,1]$ there exist $N_1>0$ and $c>0$ such that 
	\[\|U_\varepsilon\|^2_{X\otimes (S)_{ -p}} = \sum_{\gamma \in \mathcal{I}} \|u_{\gamma,\varepsilon}\|^2_{X}  \, (2\N)^{-p\gamma} <  c \, \,  \varepsilon^{-N_1} \]for  $p=\max\{p_1, p_2\}$. 
	As the  coefficients $u_{\gamma, \varepsilon}$ satisfy \eqref{estimate u gama epsilon nets} we obtain 
	\begin{equation*}
	\begin{split}\lefteqn{\|U_\varepsilon\|^2_{X\otimes (S)_{-p}}= 
		\sum_{\gamma \in \mathcal{I}} \|u_{\gamma,\varepsilon}\|^2_{C([0,T]; L^2(\R^d))} \,  (2\N)^{-p\gamma} = \sum_{\gamma \in \mathcal{I}} \left( \sup_{t\in [0,T]} \|u_{\gamma, \varepsilon} (t, \cdot)\|_{L^2(\R^d)} \right)^2\, (2\N)^{-p\gamma} }\\
	&\leq \sum_{\gamma \in \mathcal{I}} \left\{ \sup_{t\in [0,T]}  Me^{wt} \left( \|g_\gamma (\cdot)\|_{L^{2}(\R^d)} +\int_0^t \|f_\gamma (\tau,\cdot)\|_{L^2(\R^d)} \,d\tau\right) \, \varepsilon^{-N} \right\}^2 \, (2\N)^{-p\gamma}\\
	& \leq \sum_{\gamma \in \mathcal{I}} \left\{  Me^{wT} \|g_\gamma (\cdot)\|_{L^{2}(\R^d)} \, \varepsilon^{-N} + Me^{wT}\, \varepsilon^{-N} \int_0^T \|f_\gamma (\tau,\cdot)\|_{L^2(\R^d)} \,d\tau \right\}^2 \, (2\N)^{-p\gamma}\\
	& \leq \sum_{\gamma \in \mathcal{I}} \left\{ 2 M^2e^{2wT} \|g_\gamma (\cdot)\|_{L^{2}(\R^d)}^2 \, \varepsilon^{-2N} + 2M^2e^{2wT}\, \varepsilon^{-2N} \left(\int_0^T \|f_\gamma (\tau,\cdot)\|_{L^2(\R^d)} \,d\tau \right)^2\right\} \, (2\N)^{-p\gamma}\\
	&\leq 2M^2 e^{2wT}  \left\{ \sum_{\gamma \in \mathcal{I}}    \|g_{\gamma}(\cdot)\|_{L^2(\R^d)}^2   \, (2\N)^{-p\gamma}  + \sum_{\gamma \in \mathcal{I}}   \Big(\int_{0}^{T}\|f_{\gamma}(\tau, \cdot)\|_{L^2(\R^d)}\,d\tau\Big)^2 \,  (2\N)^{-p\gamma} \right\} \, \varepsilon^{-2N} \\&= 2M^2 e^{2wT}\big(I_1 + I_2) \,\, \varepsilon^{-2N}.
	\end{split}
	\end{equation*}
	Then, for  $p \geq p_2 $ we obtain
	\begin{equation*}
	\begin{split}
	I_1& = \sum_{\gamma \in \mathcal{I}}    \|g_{\gamma}(\cdot)\|_{L^2(\R^d)}^2  \, \,  (2\N)^{-p\gamma}   
	\leq \sum_{\gamma \in \mathcal{I}}    \|g_{\gamma}(\cdot)\|_{L^2(\R^d)}^2  \,(2\mathbb N)^{-p_2\gamma} 
	\end{split}
	\end{equation*}
	which is finite by \eqref{poc uslov konv}. Similarly, for  $p\geq p_1$
	\begin{equation*}
	\begin{split}
	I_2&=   \sum_{\gamma \in \mathcal{I}}   \Big(\int_{0}^{T}\|f_{\gamma}(\tau, \cdot)\|_{L^2(\R^d)}\,d\tau\Big)^2 \, (2\N)^{-p\gamma} \\
	&\leq   \sum_{\gamma \in \mathcal{I}}   \Big(\int_{0}^{T} \sup_{\tau \in [0,T]} \|f_{\gamma}(\tau, \cdot)\|_{L^2(\R^d)}\,d\tau\Big)^2 \, (2\N)^{-p\gamma} \\
	&=   \sum_{\gamma \in \mathcal{I}}   \Big(\int_{0}^{T}  \|f_{\gamma}\|_{X}\,d\tau\Big)^2 \, (2\N)^{-p\gamma} \\
	&= T^2\,    \, 
	\sum_{\gamma \in \mathcal{I}}   \|f_\gamma\|^2_X  \,\cdot  (2\N)^{-p\gamma} \\ 
	&\leq 
	\sum_{\gamma \in \mathcal{I}}   \|f_\gamma\|^2_X   (2\N)^{-p_1\gamma} 
	\end{split}
	\end{equation*}
	which is finite by \eqref{force term uslov konv}. We  conclude that for every $\varepsilon\in (0,1]$,   the stochastic process $U_\varepsilon$  is   $X\otimes (S)_{ -p}$-moderate for $p \geq \max\{p_1, p_2\}$ 
	and $U_\varepsilon$ solves  \eqref{Eq: stochastic evolution}, i.e., according to Definition \ref{moderate gsp}  $(U_\varepsilon)_{\varepsilon \in (0,1]}$ is a stochastic very weak solution to the problem \eqref{Eq: stochastic evolution}.   \hfill $\square$

\vspace{0.2cm}
\subsection{Uniqueness of the stochastic very weak solution} \label{subsec 3.2} 

As in the deterministic case the question of uniqueness for the obtained very weak solution to the stochastic parabolic  problem \eqref{Eq: stochastic evolution} can be treated in different manners.   
Let $(U_\varepsilon)_{\varepsilon \in (0,1]}$ and $(\tilde{U}_\varepsilon)_{\varepsilon \in (0,1]}$ be two very weak solutions to  the stochastic problem \eqref{Eq: stochastic evolution} that correspond to the same regularizing net $(q_\varepsilon)_{\varepsilon \in (0,1]}$  of the potential {$q\in \mathcal D'(\mathbb R^d)$}. Let $V_\varepsilon := U_\varepsilon - \tilde{U}_\varepsilon,$ $\varepsilon \in (0,1],$ i.e., $$V_\varepsilon (t,x,\omega)=\sum_{\gamma \in \mathcal{I}}v_{\gamma,\varepsilon}(t,x)H_{\gamma}(\omega):=\sum_{\gamma \in \mathcal{I}} (u_{\gamma,\varepsilon} - \tilde{u}_{\gamma,\varepsilon})(t,x) H_{\gamma}(\omega).$$ Then,  $(V_\varepsilon)_{\varepsilon \in (0,1]}$ is a very weak solution to  the stochastic homogeneous problem 	
\begin{equation*}
\label{hom stoch heat with L}
\begin{split}
\left(\frac{\partial}{\partial t}  - \mathcal L \right) V(t, x, \omega) + {q(x)} \cdot V(t, x, \omega) &= 0, 
\\
V(0, x, \omega) &=0.
\end{split}\enspace
\end{equation*}
Hence,  $(v_{\gamma,\varepsilon})_{\varepsilon \in (0,1]}$ is a very weak solution to  the corresponding homogeneous deterministic problem
\begin{equation*}
\label{hom det heat with L}
\begin{split}
\left(\frac{\partial}{\partial t}  - \L \right) v_\gamma (t, x) + {q_\varepsilon (x)} \cdot v_\gamma(t, x) &= 0, 
\\
v_\gamma(0, x) &=0, 
\end{split}\enspace
\end{equation*}
for each $\gamma\in \mathcal I$ and every $\varepsilon \in (0, 1]$. 
{From Theorem \ref{theorem estimates L infinity deterministic} and the  estimate \eqref{ee-perturbations} it follows that $v_{\gamma,\varepsilon} = 0,$ $\varepsilon \in (0,1],$ and therefore $u_{\gamma,\varepsilon} = \tilde{u}_{\gamma,\varepsilon}$, $\varepsilon \in (0,1].$ We conclude that $U_\varepsilon = \tilde{U}_\varepsilon,$ $\varepsilon \in (0,1]$, since their coefficients coincide. }
\begin{thm}
	\label{Thm: independence of regularizattion}
	Let $(q_\varepsilon)_{\varepsilon \in (0,1]}$ and $(\tilde q_\varepsilon)_{\varepsilon \in (0,1]}$ be two different regularizing nets of {$q\in \mathcal D'(\mathbb R^d)$} such that for every $\varepsilon\in (0,1]$ it holds 
	\begin{equation*} 
	\label{q negligible} 
	\| q_\varepsilon - \tilde q_\varepsilon\|_{L^{\infty}(\mathbb{R}^d)} \leq  \, C \, \varepsilon^n, \quad \forall n \in \mathbb{N}.
	\end{equation*}
	Let  $(U_{\varepsilon})_{\varepsilon \in (0,1]} \subset \mathcal{X}(L^2(\R^d)) $ and $(\tilde U_{\varepsilon})_{\varepsilon \in (0,1]} \subset \mathcal{X}(L^2(\R^d))$ be  two very weak solutions of  the stochastic  problem \eqref{Eq: stochastic evolution} which correspond to  $(q_\varepsilon)_{\varepsilon \in (0,1]}$ and  $(\tilde q_\varepsilon)_{\varepsilon \in (0,1]}$ respectively. Then, 
	for all $\varepsilon \in (0,1]$ and all $n\in \N$ there exists $c>0$ such that
	\begin{equation*}
	\|U_\varepsilon - \tilde U_\varepsilon\|_{X \otimes (S)_{-s}} \leq c\, \varepsilon^n.
	\end{equation*}
\end{thm}

{\bf Proof.} 
	Two very weak solutions $(U_{\varepsilon})_{\varepsilon \in (0,1]}$ and $(\tilde U_{\varepsilon})_{\varepsilon \in (0,1]}$  to the stochastic  initial value problem \eqref{Eq: stochastic evolution} in $\mathcal{X}(L^2(\R^d))$ that correspond to   $(q_\varepsilon)_{\varepsilon \in (0,1]}$ and  $(\tilde q_\varepsilon)_{\varepsilon \in (0,1]}$ respectively,  are $C([0,T], L^2(\mathbb R^d)) \otimes (S)_{-s}$-moderate  for some $s>1$ (see Remark \ref{Remark2}) and can be represented in the chaos expansion forms $$U_\varepsilon (t,x,\omega)= \sum_{\gamma \in \mathcal{I}} u_{\gamma,\varepsilon}(t,x) H_{\gamma}(\omega), \quad \text{and} \quad 
	\tilde{U}_\varepsilon (t,x,\omega)= \sum_{\gamma \in \mathcal{I}} \tilde{u}_{\gamma,\varepsilon}(t,x) H_{\gamma}(\omega)$$
	with 
	\[ \begin{split} \|U_\varepsilon\|_{X \otimes (S)_{-s}}  &= \sum_{\gamma \in \mathcal{I}} \|u_{\gamma,\varepsilon}\|^2_{X} (2\mathbb{N})^{-s\gamma}< \infty \quad \textnormal{and} \\
	\|\tilde{U}_\varepsilon\|_{X \otimes (S)_{-s}} &= \sum_{\gamma \in \mathcal{I}} \|\tilde{u}_{\gamma,\varepsilon}\|^2_{X} (2\mathbb{N})^{-s\gamma}< \infty.
	\end{split}\]
	Then the stochastic process $U_\varepsilon - \tilde{U}_{\varepsilon},$ $\varepsilon \in (0,1]$ has the  chaos expansion form 
	$$(U_\varepsilon - \tilde{U}_{\varepsilon})(t,x,\omega) = \sum_{\gamma \in \mathcal{I}} (u_{\gamma,\varepsilon}-\tilde{u}_{\gamma,\varepsilon})(t,x) H_{\gamma}(\omega). $$ 
	We need to show that 
	$$\|U_\varepsilon - \tilde U_\varepsilon\|^2_{X \otimes (S)_{-s}} = \sum_{\gamma \in \mathcal{I}} \|u_{\gamma,\varepsilon}-\tilde{u}_{\gamma,\varepsilon}\|^2_{X} (2\mathbb{N})^{-s\gamma}\leq c\, \varepsilon^n  \quad \textnormal{for all } n \in \mathbb{N}.$$ Since the nets $(u_{\gamma,\varepsilon})_{\varepsilon \in (0,1]}$ and $(\tilde{u}_{\gamma,\varepsilon})_{\varepsilon \in (0,1]}$ are solutions of the corresponding deterministic problems, Theorem \ref{thm:uniqueness deterministic}  implies $$\|u_{\gamma,\varepsilon}-\tilde{u}_{\gamma,\varepsilon}\|_{X} \leq   c\, \varepsilon^n  \quad \textnormal{for all } n \in \mathbb{N},$$ 
	and therefore for all $n \in \mathbb{N}$ it holds
	$$\|U_\varepsilon - \tilde U_\varepsilon\|^2_{X \otimes (S)_{-s}} = \sum_{\gamma \in \mathcal{I}} \|u_{\gamma,\varepsilon}-\tilde{u}_{\gamma,\varepsilon}\|^2_{X} (2\mathbb{N})^{-s\gamma}\leq  \sum_{\gamma \in \mathcal{I}} c \varepsilon^n (2\mathbb{N})^{-s\gamma}   = C \,\, \varepsilon^n ,$$ 
	with $C= \sum_{\gamma \in \mathcal{I}} c \, (2\mathbb{N})^{-s\gamma} < \infty$  since $s>1.$   \hfill $\square$

\vspace{0.2cm}
\subsection{Consistency} \label{subsec 3.3}
In this section we are interested in consistency of the stochastic very weak solution obtained in Theorem \ref{Thm 1: exist} with the stochastic (classical weak)   solution obtained in Theorem \ref{stochastic weak sol} in the following sense. 
\begin{thm}
	\label{Thm 2: consistency} 
	Consider the stochastic parabolic  problem \eqref{Eq: stochastic evolution}, let the potential $q\in L^\infty(\mathbb R^d)\cap C(\mathbb R^d)$ and  
	let the assumptions $(H1)$-$(H3)$ hold. Let $V\in \mathcal{X}(L^2(\R^d))$ be the solution  to the stochastic parabolic  problem obtained in Theorem \ref{stochastic weak sol} and let $(U_\varepsilon)_{\varepsilon\in (0,1]} \subset \mathcal{X}(L^2(\R^d))$  be the very weak solution to the stochastic  problem obtained in Theorem \ref{Thm 1: exist}.  Then, for some $s>1$
	\begin{equation}
	\label{consistency conv}
	\|U_\varepsilon - V \|_{X \otimes (S)_{ -s}} \to 0 \qquad \text{as} \qquad \varepsilon \to 0.
	\end{equation} 
\end{thm}

{\bf Proof.} 
	Since it is assumed that $(U_\varepsilon)_{\varepsilon\in (0,1]}$  is a very weak solution to the  problem \eqref{Eq: stochastic evolution},  then for each $\varepsilon\in (0,1]$ it is a stochastic process in $\mathcal{X}(L^2(\R^d))$ represented in the form 
	\begin{equation*}
	U_{\varepsilon}(t, x, \omega) = \sum_{\gamma \in \mathcal I} u_{\gamma, \varepsilon} (t, x) \,  H_\gamma(\omega). 
	\end{equation*}
	Moreover, $(U_\varepsilon)_{\varepsilon\in (0,1]}$ is $X \otimes (S)_{ -p}$-moderate, which is equivalent to the net of coefficients $(u_{\gamma, \varepsilon})_{\varepsilon \in (0,1]}$ being $X$-moderate for each $\gamma \in \mathcal I$. Recall, the coefficients $u_{\gamma, \varepsilon}$ for all $\gamma\in \mathcal I$ and all $\varepsilon\in (0,1]$ solve 
	\begin{equation}
	\label{reg det pde}
	\begin{split}
	\left( \frac{\partial}{\partial t}  -  {\L} \right) u_\gamma(t, x) + {q_\varepsilon (x)} \cdot  u_\gamma(t, x) & = f_\gamma(t, x ), \quad t\in (0,T], \\
	u_\gamma(0, x) &= g_\gamma(x),
	\end{split}
	\end{equation}  where   $q_\varepsilon=q*\varphi_\varepsilon \in {C}^\infty(\mathbb R^d)$.  
	On the other side, the process $V\in \mathcal{X}(L^2(\R^d))$ is a  classical weak solution to \eqref{Eq: stochastic evolution}.
	It has  chaos expansion representation of the form 
	\begin{equation*}
	\label{process V chaos exp}
	V(t, x, \omega) =\sum_{\gamma \in \mathcal I} v_{\gamma} (t, x) H_\gamma(\omega), 
	\end{equation*}
	where its coefficients $v_\gamma$ for each $\gamma\in \mathcal I$ are (classical) weak solutions to the deterministic problem 
	\begin{equation}
	\label{coefficients of V}
	\begin{split}
	\frac{\partial}{\partial t} v_\gamma(t, x)  -  {\L}   v_\gamma(t, x) + {q (x)} \cdot  v_\gamma(t, x) & = f_\gamma(t, x ), \quad t\in (0,T], \,   x\in  \mathbb R^d,   
	\\
	v_\gamma(0, x) = g_\gamma(x),
	\end{split}
	\end{equation} with $q\in L^\infty(\mathbb R^d)\cap C(\mathbb R^d)$,  
	such that  for $r\geq 0$  it holds 
	\[\sum_{\gamma \in \mathcal I} \|v_{\gamma}\|^2_X \, (2\mathbb N)^{-r\gamma } < \infty  .\]
	In order to prove \eqref{consistency conv} we have to prove 
	\[\|U_\varepsilon - V \|_{X \otimes (S)_{ -s}} =\sum_{\gamma\in \mathcal I} \|u_{\gamma, \varepsilon} - v_\gamma\|^2_X \, (2\mathbb N)^{-s\gamma} \to 0 \quad \text{as} \quad \varepsilon \to 0.\]
	This is true since by Theorem \ref{thm:consistency deterministic} the very weak solution   $u_{\gamma,\varepsilon}$ to  \eqref{reg det pde}  is consistent to the weak solution  $v_\gamma$ to \eqref{coefficients of V} for each $\gamma\in \mathcal I$, i.e.,   it follows that 
	$\|u_{\gamma, \varepsilon}- v_\gamma\|_X \to 0$ as $\varepsilon \to 0$, for all $\gamma\in \mathcal I$. 	
	\hfill $\square$ 

\vspace{0.2cm}
\section{An example, conclusions and further extensions} \label{sec5}
To illustrate the proposed method  we consider an example of a stochastic heat equation with singular potential.  
Some notations and basic concepts introduced earlier that are not much used throughout the manuscript will be used here and we refer reader to the Subsection \ref{NotNot} for recalling the details.       

In the problem \eqref{Eq: stochastic evolution} let  the operator $\L$ be the Laplace operator over the domain $D=H^1_0(\R^d)$. Then  $\frac{\partial}{\partial t} - \Delta$ may be seen as the semigroup evolving law of a rescaled Brownian motion. Let further the potential $q$ be Dirac delta distribution in space, let  the force term $F$ be the time-white noise process  $W\in C^k([0,T]) \otimes (S)_{-1}$ defined in Example  \ref{Example0} (iv), and the initial condition $G$ be a non-zero Gaussian generalized stochastic process in $H_0^1(\R^d) \otimes (S)_{-1}$. Thus, for such data and $t\in (0,T]$, $x\in \mathbb \R^d$, $ \omega\in \Omega$, we consider the Cauchy problem
\begin{equation}
\label{Example1}
\begin{split}
\left(\frac{\partial}{\partial t}  - \Delta \right) \, U(t, x, \omega) + \delta(x) \cdot U(t, x, \omega) &= W(t,\omega), \enspace \\
U(0, x, \omega) &= G(x,\omega).
\end{split}
\end{equation}
The white noise process  $W$ has the chaos expansion representation given by \eqref{whitenoise} while 
the Gaussian generalized stochastic process $G$  has the chaos expansion representation (see \cite{LS2017})
$$G(x, \omega) = g_0 (x) + \sum_{k=1}^\infty g_{e(k)}(x) \, H_{e(k)}(\omega),$$ 
where $ g_0= \mathbb{E}( G)$ is the expectation of $G$, $e(k) $ denotes $k$th unit vector being the sequence of zeros and having $1$ as the $k$th component, and $H_{e(k)}(\omega)$ are Fourier-Hermite polynomials given by \eqref{Halpha}.

We assume that the solution to \eqref{Example1}  is given by \eqref{sol U}, more precisely in the form
\begin{equation*}
\label{oblik resenja koji trazimo1}
U(t, x, \omega) = u_0 (t,x) + \sum_{k=1}^{\infty} u_{e(k)}(t, x) \, H_{e(k)}(\omega) + \sum_{|\gamma|>1} u_\gamma(t,x) H_\gamma(\omega), 
\end{equation*} where $u_0$ is the expectation of  $U$, i.e. $u_0=\mathbb{E}(U)$.
The chaos expansion method implies that for $|\gamma|=0$  the expectation $u_0$  satisfies the deterministic PDE
\begin{equation}
\label{Example1: the system of pde nulti}
\begin{split}
\frac{\partial}{\partial t} u_{0}(t, x) - \Delta u_{0}(t, x) + \delta(x) \,   u_{0}(t, x)&= 0, \enspace \\
u_{0}(0, x)&=g_{0}(x).
\end{split}
\end{equation}
For $|\gamma|=1$,   $\gamma = e(k) $,  $k \in \mathbb{N}$,  the coefficients $u_{e(k)}$ satisfy 
\begin{equation}
\label{Example1: the system of pde}
\begin{split}
\frac{\partial}{\partial t} u_{e(k)}(t, x) - \Delta u_{e(k)}(t, x) + \delta(x) \,   u_{e(k)}(t, x)&= \xi_k(t), \enspace \\
u_{e(k)}(0, x)&=g_{e(k)}(x).
\end{split}
\end{equation}
For $|\gamma|>1$ the coefficients  $u_\gamma$ satisfy the homogeneous problem of the form \eqref{Example1: the system of pde nulti} with $g_\gamma=0$. 
The regularizing net of the Dirac delta  distribution   is any mollifying net of smooth functions $\delta_\varepsilon$, for example  
$$\delta_\varepsilon (x)  = \frac{1}{\varepsilon^d} \varphi \left(\frac{x}{\varepsilon} \right),$$
with    $\varphi$ being a mollifier (see Subsection \ref{SubSecReg}). 
The sequence of regularized problems corresponding to $|\gamma|=0$  and to problem \eqref{Example1: the system of pde nulti} is given by 
\begin{equation}
\label{Example1: regularized pde nulti}
\begin{split}
\frac{\partial}{\partial t} u_{0}(t, x) - \Delta u_{0}(t, x) + \delta_\varepsilon(x)  u_{0}(t, x)&= 0, \enspace \\
u_{0}(0, x)&=g_{0}(x),
\end{split}
\end{equation} 
while  the sequence of regularized problems corresponding to $|\gamma|=1$ and to problems \eqref{Example1: the system of pde} is given by 
\begin{equation}
\label{Example1: regularized pde}
\begin{split}
\frac{\partial}{\partial t} u_{e(k)}(t, x) - \Delta u_{e(k)}(t, x) + \delta_\varepsilon (x)  u_{e(k)}(t, x)&= \xi_k(t), \enspace \\
u_{e(k)}(0, x)&=g_{e(k)}(x).
\end{split}
\end{equation}
The solution to the homogenous regularized problem
\eqref{Example1: regularized pde nulti} is given by 
$$u_{0,\varepsilon} (t,x)=   S_t g_{0} (x), \quad t\in [0,T], \, x \in \mathbb{R}^d,
$$
where $S_t$ is the semigroup generated by perturbed Laplace operator $\Delta - \delta_\varepsilon (x) Id$ and it is given by
\begin{equation}\label{Sodt}
S_t= e^{t \delta_\varepsilon (x)} T_t, \quad t \geq 0,
\end{equation}
with $\{T_t\}_{t\geq 0}$ being $C_0$-semigroup generated by the  Laplace operator, see \cite{Pazy1983}.
The solution to the regularized problem
\eqref{Example1: regularized pde} is given by
$$u_{e(k),\varepsilon} (t,x)=   S_t g_{e(k)} (x)+\int_0^t S_{t-s} \xi_k(s)\,ds, \quad t\in [0,T], \, x \in \mathbb{R}^d,
$$ where again,  $S_t$ is  the semigroup given by \eqref{Sodt}.  For $|\gamma|>1$  the coefficients $u_\gamma$ satisfy the regularized homogeneous problem  \eqref{Example1: regularized pde nulti} with $g_\gamma=0$, and therefore  for $|\gamma|>1$ coefficients $u_{\gamma,\varepsilon}=0$. 
Thus, solution to \eqref{Example1} is given by
\begin{equation*} 
U(t, x, \omega) = S_t g_{0} (x) + \sum_{k=1}^{\infty} \left(S_t g_{e(k)} (x)+\int_0^t S_{t-s} \xi_k(s)\,ds\right) \, H_{e(k)}(\omega).
\end{equation*}

\vspace{0.3cm}
\subsection{Conclusions and further extensions}
This paper brings the following  novelties.  
The notion of stochastic very weak solutions  for the given stochastic parabolic  initial value problem \eqref{Eq: stochastic evolution} with irregular potential is introduced. In Theorem \ref{Thm 1: exist} and Theorem \ref{Thm: independence of regularizattion} we proved the existence and uniqueness of a  stochastic very weak solution to \eqref{Eq: stochastic evolution},  while in Theorem \ref{Thm 2: consistency} we proved that  the stochastic very weak solution is consistent with a  stochastic (classical) weak solution.  In addition,  Theorem \ref{thm extended Rozha}  where we have proved the existence and uniqueness of a very weak solution to the  deterministic parabolic  equation with singular potential, extends  the  results from \cite{ARST2020c}, where the authors considered the  heat equation with singular potential without a force term.

As the proposed method   works very well in this simplest case, we aim to apply it to the general problem (\ref{Eq: general}).
For $q$ being time dependent and bounded function, the theory of semigroups  guarantees the existence of an evolution system $E(s,t)$ in the place of the semigroup $T(t)$. The question of what happens if $q$ is irregular in time arises. A possible way to overcome this difficulty would be to apply the regularization procedure from \cite{NPR2005,RO1999}. In our further work, we are going to  take the advantage of   the  very weak solution approach developed in this paper. On the other hand, stochastic irregularities will be considered in white noise analysis setting where the product is interpreted as the Wick product. Also we are going to apply our method   to semilinear stochastic problems.  Finally, as the concept of very weak solutions is well adapted in numerical studies of different classes of PDEs with singularities \cite{ARST2021,ARST2020b,ARST2020c,MRT2019a,MRT2019b},  we aim to adapt the proposed method  for numerical approximations of   different physical phenomena modelled by  \eqref{Eq: general} as it was done for example in \cite{LM2017}. 

To conclude, we single out the highlights of the paper. The very weak solution for the stochastic parabolic initial value problem with irregular potential is introduced. For the analysis, the chaos expansion method from the white noise analysis and the concept of very weak solutions to treat singularities are combined. The existence and uniqueness of a stochastic very weak solution are proved, and the consistency of obtained solution with classical notions of the solutions is shown. The newly proposed method could be  applied  to a wide classes of stochastic partial differential equations with singularities and it is suitable for their numerical analysis.


\end{document}